\documentclass[twoside]{bourlarge}
\usepackage{mazliak_h}

\usepackage{amssymb}
\usepackage{amsmath}
\usepackage{graphicx}

\usepackage[dvips]{epsfig}

\usepackage{color}

\newcommand{\ind}{1{\hspace{-0.1cm}}{\rm{I}}}

\newcommand{\cqfd} {\hbox {\unskip \kern 6pt \penalty 500\raise -2pt \hbox
{\vrule \vbox to 5pt {\hrule width 4pt \vfill \hrule }\vrule
}\par }}

\newlength\jataille
\newcommand{\figgauche}[3]%
{
\jataille=\textwidth\advance\jataille by -#1
\advance\jataille by -.5cm
\begin{minipage}[a]{#1}
\includegraphics[width=#1]{#2}
\end{minipage}
\vskip2mm
\begin{minipage}[a]{\jataille}
\footnotesize #3 \normalsize
\end{minipage}
}

\usepackage{amsfonts}

\font\tencyr=wncyr10          
\def\cyr{\tencyr

}

\def \ikr{\accent"24 i}

\def\ppref
#1
{$\bullet $
\ref{#1}
}
\newcommand{\lettre}[1]{
\refstepcounter{section}
\vskip 5mm
\centerline{Lettre \thesection}
\sectionmark{Lettre \thesection}
\vskip 3mm }

\usepackage{ulem}
\usepackage[english]{babel}
\usepackage[cp1250]{inputenc}
\makeindex

\begin{document}
\Large

\title{Poincar{\'e}'s Odds}

\author{Laurent {\sc Mazliak} \\
Laboratoire de Probabilit{\'e}s et Mod{\`e}les Al{\'e}atoires \\ Universit{\'e} Pierre et Marie Curie\\ Paris, France }

\maketitle

\begin{abstract}
This paper is devoted to Poincar{\'e}'s work in probability. Though the subject does not represent a large part of the mathematician's achievements, it provides significant insight into the evolution of Poincar{\'e}'s thought on several important matters such as the changes in physics implied by statistical mechanics and  molecular theories. After having drawn the general historical context of this evolution, I focus on several important steps in Poincar{\'e}'s texts dealing with probability theory, and eventually consider how his legacy was developed by the next generation.

\end{abstract}

\section*{Introduction}

In 1906, Poincar{\'e} signed one of the most unusual texts of his
scientific career \cite{AppelDarbouxPoincare1906}, a report written
for the Cour de Cassation in order to eventually close the Dreyfus
case. In 1904, ten years after the arrest and the degradation of the unfortunate captain, the French government
decided to bring an end to this lamentable story which had torn French
society apart for years, and obtained the rehabilitation of the young
officer who had been so unjustly martyred. As is well known, Dreyfus
was convicted of treason in 1894 in the absence of any material proof,
during a rushed and unbalanced trial orchestrated by the military
hierarchy, whose aim was to identify a culprit as soon as
possible. The only concrete document submitted as evidence by the
prosecution was the famous {\it bordereau}, found in a wastepaper
basket in the German Embassy in Paris and briefly scrutinized by
several more or less competent experts. Among them, Alphonse Bertillon
played a specially sinister role and became Dreyfus' most obstinate
accuser. He built a bizarre, scientific-sounding theory of
self-forgery of the {\it bordereau}, purporting to prove Dreyfus'
guilt. Bertillon became trapped by his own theoretical construction,
more because of overweening self-confidence and stupidity than by a
real partisan spirit. When the Affair broke out at the end of the
1890s, and the political plot became obvious and Dreyfus' innocence
apparent, Bertillon continued to elaborate his theory, rendering it
more and more complex. This frenzy of elaboration brought him an
avalanche of trouble, and nearly put an end to his career in the Paris
police department. However, for the trial in the {\it Cour de
  Cassation}, in order to silence the last dissenting voices which
might be raised, it was decided not to ignore Bertillon's rantings,
but to ask incontestable academic authorities to give their opinion
about the possible value of the self-proclaimed expert's
conclusions. Three mathematicians were called for that purpose, Paul
Appell, Gaston Darboux and Henri Poincar{\'e}, who jointly signed the
report for the Cour de Cassation in 1906. Nevertheless, it was well
known that only the latter worked on the document, a tedious task
Poincar{\'e} that undertook with integrity, but somewhat
grudgingly. Besides, this was not Poincar{\'e}'s first intervention in the
Dreyfus case: during the judicial review of 1899 in Rennes of the
decision of 1894, Painlev{\'e} read a letter of Poincar{\'e} criticizing
Bertillon's work for its lack of scientific foundation.

This story is well known and has been narrated in detail many times
(\cite{MansuyMazliak2005}, \cite{MansuyMazliak2011},
\cite{Rollet1999}) and I shall not dwell on it further. However, one
may ask why it was Poincar{\'e} who was called upon to carry out this
task. The first answer that comes to mind is the following: in 1906, Poincar{\'e}, who
was then fifty-two years old, was without a doubt the most prominent among French
scientists. He had, moreover, another characteristic: his name was
familiar to a sizeable audience outside the scientific community. It
was well known that he had been awarded the important prize of the Swedish king
Oscar II in 1889; the publication of popular books on the
philosophy of science increased his name recognition, and he was also
famous due to several essays published in newspapers, for
instance during the International Congress of Mathematicians in Paris
in 1900, over which Poincar{\'e} presided. To seek out such a
person in order to crush the insignificant but noisy Bertillon was
therefore a logical calculation on the government's part. However, a
second and more hidden reason probably played a part as well. The report
for the {\it Cour de Cassation} \cite{AppelDarbouxPoincare1906} opens
with a chapter whose title was quite original in the judicial
litterature: {\it Notions on the probability of causes}; it contains a
brief exposition of the principles of the Bayesian method. Bertillon
had indeed pretended to build his system on the methods and results of
the theory of probability, and answering him was only
possible by confronting the so-called expert with his own weapons. It
was therefore necessary not only to seek out a scientific celebrity for the
job, but also someone whose authority in these matters could not be
challenged. In 1906, Poincar{\'e} was regarded by everyone
as the leading specialist in the mathematics of randomness in
France. While by then he was no longer the Sorbonne Chair of
Calculus of Probability and Mathematical Physics, he had held the Chair
for some ten years (it was in fact the first position he obtained in
Paris), had accomplished an impressive amount of work in mathematical
physics during this period (we shall return at length to this subject
later), and had published in 1896, shortly before leaving the chair, a
treatise on the Calculus of Probability which, in 1906, was still the
main textbook on the subject in French. Moreover, several texts had
presented his thoughts on the presence of randomness in modern physics
to a large cultivated audience, in particular his {\it Science and
  Hypothesis} (\cite{Poincare1902}), which went through several
editions. It
was therefore as a specialist in the calculus of probability that
Poincar{\'e} was sought out by the judicial authorities, in order to
dispose of the Dreyfus case once and for all.

If we go back in time fifteen years before this event, we notice a
striking contrast. Since 1886, although Poincar{\'e} held the
aforementioned Sorbonne chair, he had probably mainly seen
in the chair's title the words {\it Mathematical Physics}. In 1892, he
published an important textbook on Thermodynamics \cite{Poincare1892}
based on the lectures he had read at the Sorbonne several years
before. A Poincar{\'e} publication would of course not go unnoticed, and
one attentive reader had been the British physicist Peter Guthrie Tait
(1831-1901). Tait had been very close to Maxwell and was one of the
most enthusiastic followers of his work. He wrote a review of
Poincar{\'e}'s book for the journal {\it Nature} (\cite{Tait1892}); the
review was quite negative despite the obvious talent Tait recognized
in his young French colleague. 

Poincar{\'e}, wrote Tait, introduced beautiful
and complex mathematical theories in his textbook but often to the
detriment of the physical meaning of the situations he studied. The
most important reproach of the British physicist was the fact that
Poincar{\'e} remained absolutely silent about the statistical theories of
thermodynamics, leaving in the shadows the works of Tait's friend and
master Maxwell. Tait wrote:

\begin{quote}`But the most unsatisfactory part of the whole work is,
  it seems to us, the entire ignoration of the true (i.e. the
  statistical) basis of the Second Law of Thermodynamics. According to
  Clerk-Maxwell (Nature, xvii. 278) ``The touch-stone of a treatise on
  Thermodynamics is what is called the Second Law.''  We need not
  quote the very clear statement which follows this, as it is probably
  accessible to all our readers. It certainly has not much resemblance
  to what will be found on the point in M. Poincar{\'e}'s work: so little,
  indeed, that if we were to judge by these two writings alone it
  would appear that, with the exception of the portion treated in the
  recent investigations of v. Helmholtz, the science had been
  retrograding, certainly not advancing, for the last twenty
  years.'\end{quote}

Poincar{\'e} wrote an answer and sent it to {\it Nature} on 24 February
1892.  It was followed, during the first semester of 1892, by six
other letters between Tait and Poincar{\'e}\footnote{Accessible freely on the website of the Archives Poincar\'e in Nancy at the address 
http://www.univ-nancy2.fr/poincare/chp/hpcoalpha.xml }; they are rather sharp, each
sticking to his position. On March 17th, Poincar{\'e} wrote the following
comment about the major criticism made by Tait:

\begin{quote} `I left completely aside a mechanical explanation of the
  principle of Clausius that M. Tait calls ``the true
  (i.e. statistical) basis of the Second Law of Thermodynamics.'' I
  did not speak about this explanation, which besides seems to me
  rather unsatisfactory, because I desired to remain absolutely
  outside all the molecular hypotheses, as ingenious they may be; and
  in particular I said nothing about the kinetic theory of
  gases.'\footnote{ `J'ai laiss{\'e} compl{\`e}tement de c{\^o}t{\'e} une explication
    m{\'e}canique du principe de Clausius que M. Tait appelle ``the true
    (i.e. statistical) basis of the Second Law of Thermodynamics.'' Je
    n'ai pas parl{\'e} de cette explication, qui me para{\^\i}t d'ailleurs
    assez peu satisfaisante, parce que je d{\'e}sirais rester compl{\`e}tement
    en dehors de toutes les hypoth{\`e}ses mol{\'e}culaires quelque
    ing{\'e}nieuses qu'elles puissent {\^e}tre; et en particulier j'ai pass{\'e}
    sous silence la th{\'e}orie cin{\'e}tique des gaz'.}\end{quote}

One thus observes that Poincar{\'e}, in 1892, had a very negative view of
statistical mechanics, the very science of matter in which
probabilities played the greatest role. However, Poincar{\'e} would not
have been Poincar{\'e} if, once shown a difficulty, he did not take the
bull by the horns and to try to tame it. A first decision was the
resolution, during the next academic year, 1893, to lecture on the
kinetic theory of gases. And indeed, in the Sorbonne syllabus for that
year, we see that Poincar{\'e} had transformed his lectures into {\it
  Thermodynamics and the kinetic theory of gases}.

In 1894 Poincar{\'e}'s first published paper on the kinetic theory of
gases \cite{Poincare1894}, to which we shall return later. If one
still observes a good deal of skepticism in it, or at least
reservation about these novelties in randomness, a change of attitude
was taking place. This same academic year 1893-1894, Poincar{\'e}
eventually decided for the first time to teach a course on the
calculus of probability at the Sorbonne, which subsequently became the
basis of his book of 1896, the same year when he exchanged his chair
for that of celestial mechanics. Moreover, in the following
years, reflections on the mathematics of randomness became more
frequent in his writings, up to the publication of his more
philosophical works, which acknowledged the integration of the theory
of probability among Poincar{\'e}'s mathematical tools. By 1906, as
already noted, the transition had been completed, especially as new
elements, such as Einstein's just published theory of Brownian motion,
made even more necessary the increasing presence of probability in
scientific theories.
 
The present paper concerns the probabilistic aspects of Poincar{\'e}'s
enormous production, aspects of which remain somewhat limited in size. A
non-negligible challenge in dealing with this subject is that
probability penetrated Poincar{\'e}'s work almost by force, forcing his
hand several times, his main achievement consisting in building dykes
so that the mathematician might venture with a dry foot on these
rather soft lands. We shall see later that his successor Borel had a
somewhat different attitude toward choosing to apply probability
theory in many domains, reflecting the fact that Borel had encountered
probability in a more spectacular way than Poincar{\'e}.
 
Despite this limited contribution, Poincar{\'e} succeeded in leaving a
significant heritage which would later prove important.  Above all,
his most decisive influence may have been to allow probability theory
to regain its prestige in France again, after occupying a rather miserable
position in the French academic world for more than half a century.
 
The aim of the present paper is therefore threefold, and this is
reflected in its three parts. In the first part, we focus on
Poincar{\'e}'s evolution in the fifteen years we have just pointed to, in
order to define how this progressive taming of the mathematics of
randomness occurred. In the second section, I examine in detail
several of Poincar{\'e}'s works involving probability considerations in
order to give an idea not only of Poincar{\'e}'s style, but also of the
intellectual basis for the evolution in his thinking. And finally, in
the last section, I comment on the persons who recovered the heritage
of Poincar{\'e} in this area, and on how they extended it.
 
Many of the topics discussed in this paper have been covered
previously. In particular, I shall refer several times to Sheynin's
text (\cite{Sheynin1991}), to von Plato's publications
(\cite{Vonplato1991}, \cite{Vonplato1994}) as well as to several
papers by Bru, some published and some not (\cite{Bru2003},
\cite{Bru2013}). Let us me also mention the nice book \cite{Gargani2012} where the author deals with 
Poincar\'e's works on probability with an epistemological approach.

 \section*{Acknowledgement} 
 I warmly thank the organizers of the S{\'e}minaire Poincar{\'e} on 24
 November 2012 in Paris, for the opportunity to end the Poincar{\'e} year
 2012 with a closer reading of his work on probability, work which
 until then I knew only superficially. I hope to have succeeded in the
 paper to present my understanding of a great mind's intellectual
 peregrination, which led him to such a considerable evolution in his
 thought on the subject. Special thanks go to Bertrand Duplantier who
 made many comments on a first version of the paper. Next, I thank
 Xavier Prudent who helped me to read the 1894 paper on the kinetic
 theory of gases; needless to say, the remaining obscurities on the
 topic are entirely mine.  It is also a pleasure to express all my
 gratitude to Sandy Zabell and Scott Walter for their generous help in editing the
 English of the present paper; here also I take the complete
 responsibility for all the infelicities which may remain in the text.{
   After the presentation of this paper at the Institut Henri Poincar{\'e}
   on 24 November 2012, I had an interesting echange with Persi
   Diaconis; I hope he will find the paper improved by his comments.}
 And finally, I feel it necessary to say how much Bernard Bru would
 have deserved to write this chapter. I dedicate this work with
 gratitude to the memory of an honest man, our master and friend Marc Barbut, whom we miss
 so much since a cold morning of December 2011\dots

\section{First part: the discovery of probability}

In the beginning was the Chair. For, as we will see, it was
not only by a fortuitous scientific interest that Poincar{\'e} first came
to be involved in probability theory, but there was also a very
specific academic situation in which Hermite, the major mathematical
authority in 1880s France, crusaded for the career of his three {\it
  mathematical stars}: Paul Appell (his nephew), {\'E}mile Picard (his
son-in-law) and of course Henri Poincar{\'e}, the latter considered by him
to be the most brilliant, though he did not belong to his family (to
the great displeasure, so he wrote to Mittag-Leffler, of Madame
Hermite). A remarkable change took place at the Sorbonne at this time
when, in just a few years, almost every holder of the existing chairs
in mathematics and physics died: Liouville and Briot in 1882, Puiseux
in 1883, Bouquet and Desains in 1885, Jamin in 1886, leaving an open
field for a spectacular change in the professorships. At the end of
the whirlwind, the average age of the professors of mathematics and
physics at the Sorbonne had been reduced by eighteen years! It was
therefore clear for everyone, and first of all Hermite, that there was
a risk that the situation might become fixed for a long time, and that
it was therefore necessary to act swiftly and resolutely in favor of
his prot{\'e}g{\'e}s. The three of them were, indeed, appointed in Paris
during these years. This episode is narrated in detail in
\cite{Atten1988} and I here touch only upon the most important aspects
insofar as they relate to our story.

The Chair of Calculus of Probability and Mathematical Physics,
occupied until 1882 by Briot, had been created some thirty years
before, after numerous unsuccessful attempts by Poisson, but in a form
different than the one pursued by the latter. The joining of
mathematical physics to probability had been decided in order to
temper the bad reputation of the theory of probability in the 1840s in
France, notably due to some of the work of Laplace and, above all, of
Poisson himself dealing with the application of probability in the
judicial domain. The book of Poisson \cite{Poisson1837} ignited a
dispute in the Academy in 1836, when Dupin and Poinsot harshly
contested Poisson's conclusions, and the philosophers led by Victor
Cousin made a scene {\it in the name of the sacred rights of liberty
  against the claims advanced by the mathematicians in order to
  explain how events occur in social issues}. John Stuart Mill summed
up the situation by qualifying the application of probability to judicial
problems as the {\it scandal of mathematics}. This polemical
contretemps tarred the calculus of probability in France, leaving it
with a most dubious reputation.

Although Lippmann was nominated in 1885 for the chair, Jamin's
unexpected demise freed the chair of Research in physics of which
Lippmann took hold immediately, leaving a place for Poincar{\'e}. One
might initially be surprised by this choice. First, because Poincar{\'e}
was preferred to other candidates, of course less gifted, but whose
area of research was closer to that of the position. One would rather
have expected to see Poincar{\'e} in one of the chairs in Analysis, for
instance. In 1882, when Briot died and the great game began, Poincar{\'e}
had no published work in physics (let alone the calculus of
probability). And yet, among the candidates who were eliminated, there
was for instance Boussinesq who had to his credit significant
contributions in that domain, and who openly stated his intention of
revitalizing the teaching of probability which seemed to him to be in
a poor state. Not without irony, Boussinesq later became Poincar{\'e}'s
successor in this same chair when the latter changed over to the Chair
of Celestial Mechanics.

Thus Poincar{\'e} was nominated in 1886, without any real title for the
position, and one may think that it had really been chance which led
him there. However, as shown in Atten's fine analysis
\cite{Atten1988}, a number of indications demonstrate that Poincar{\'e} had
indeed really desired this particular position. Examining his lectures
in 1887-1888 shows that he already had a profound knowledge of
physical theories and, moreover, several passages in his
correspondence at this time show his sustained interest in
questions of physics. This shows that his attitude had not been purely
opportunistic and that the chair pleased him. Also, Hermite, apart
from the desire of supporting his prot{\'e}g{\'e}, seems to have made a
thoughtful bet in nominating him for this somewhat unexpected
position. Knowing Poincar{\'e}'s acute mind, it was not unreasonable to
anticipate spectacular achievements by him in this position. What
followed, as we know, bore out Hermite's expectation\dots

The theory of probability, as I said, did not seem to concern the new
professor at the beginning of his tenure, and he taught courses on
several different physical theories. In 1892, he published his
lectures on thermodynamics, delivered in 1888-1889 \cite{Poincare1892},
producing a book, as noted earlier, that was sharply criticized by
Tait. Poincar{\'e} therefore decided to look into the questions raised by
the kinetic theory of gases, in particular because he had just read a
communication by Lord Kelvin to the Royal Society containing several
fundamental criticisms on Maxwell's theory \cite{Kelvin1892}. Perhaps
Poincar{\'e} had been especially eager to read this paper because it might
provide a powerful argument in his controversy with Tait. However, the
affair took another direction, revealing the mathematician's profound
scientific honesty. At the beginning of the paper \cite{Poincare1894},
published in 1894, though he again expressed some skepticism,
Poincar{\'e}, whose conventionalism was formed during these years, seemed
already half convinced of the possible fecundity of Maxwell's theory.

\begin{quote}`Does this theory deserve the efforts the English devoted
  to it? One may sometimes ask the question; I doubt that, right now,
  it may explain all the facts we know. But the question is not to
  know if it is true; this word does not have any meaning when this
  kind of theory is concerned. The question is to know whether its
  fecundity is exhausted or if it can still help to make
  discoveries. And admittedly we cannot forget that it was useful to
  M.~Crookes in his research on radiant matter and also to the
  inventors of the osmotic pressure. One can therefore still make use
  of the kinetic hypothesis, as long as one is not fooled by
  it.'\footnote{ `Cette th{\'e}orie m{\'e}rite-t-elle les efforts que les
    Anglais y ont consacr{\'e}s ? On peut quelquefois se le demander; je
    doute que, d{\`e}s {\`a} pr{\'e}sent, elle puisse rendre compte de tous les
    faits connus. Mais il ne s'agit pas de savoir si elle est vraie;
    ce mot en ce qui concerne une th{\'e}orie de ce genre n'a aucun
    sens. Il s'agit de savoir si sa f{\'e}condit{\'e} est {\'e}puis{\'e}e ou si elle
    peut encore aider {\`a} faire des d{\'e}couvertes. Or, on ne saurait
    oublier qu'elle a {\'e}t{\'e} utile {\`a} M.~Crookes dans ses travaux sur la
    mati{\`e}re radiante ainsi qu'aux inventeurs de la pression
    osmotique. On peut donc encore se servir de l'hypoth{\`e}se cin{\'e}tique,
    pourvu qu'on n'en soit pas dupe.' (\cite{Poincare1894},
    p.513)}\end{quote}

We shall come back in the second section to this article of 1894 and
to the new formulation of the ergodic principle it contains, in which
Poincar{\'e} introduced the restriction of exceptional initial states. Let
us only mention here the following point: Poincar{\'e} seemed to have
found this idea in his previous work on the three body problem for
which he had obtained the prize of the King of Sweden in 1889. In the
memoir presented for the prize, he had indeed proved a recurrence
theorem concerning the existence of trajectories such that the system
comes back an infinite number of times in any region of the space, no
matter how small. The next year, for his paper in {\it Acta
  Mathematica}, Poincar{\'e} added a probabilistic extension of his
theorem where he showed that the set of initial conditions for which
the trajectories come back only a finite number of times in the
selected region has probability zero (\cite{Poincare1890}, p.71-72);
this passage was included some years later in his treatise of new
methods for celestial mechanics \cite{Poincare1892b} in a section
simply called `Probabilities'.

During the academic year 1893-1894, Poincar{\'e} prepared his first course
of probability for his students of mathematical physics
(\cite{Poincare1896}), published in 1896 and transcribed by Albert
Quiquet, a former student of the {\'E}cole Normale Sup{\'e}rieure who entered
the institution in 1883, before becoming an actuary, and who had
probably been an attentive listener of the master's voice. In the
manner of the standard probability textbooks of the time, it does not
present a unified theoretical body but rather a series of questions
that Poincar{\'e} tried to answer (the main ones, which occupy the bulk of
the volume, concern the theory of errors of measurement - a subject we shall
return to in the next section). The book, in its edition of
1896, is a natural successor to Bertrand's textbook
\cite{Bertrand1888}, which up to then was the standard textbook (see
\cite{Bru2006}). In comparison to Bertrand's book, Poincar{\'e}'s
consolidates the material in several interesting ways, and this aspect
is even more obvious in the second edition \cite{Poincare1912},
completed by Poincar{\'e} some months before his death in July, 1912.

The mathematician seemed now convinced that it was no longer possible
to get rid of probability altogether in science, so he decided instead
to make the theory as acceptable as possible to the
scientist. Poincar{\'e} decided to devote considerable effort towards that
end, especially by writing several texts lying half-way between
popularization (with the meaning of writing a description of several
modern concepts using as little technical jargon as possible) and
innovation. Two texts are of particular importance: that of 1899
(\cite{Poincare1899}) - reprinted as a chapter of \cite{Poincare1902}
- and that of 1907 (\cite{Poincare1907}) - reprinted again as
\cite{Poincare1908} and then as a preface to the second edition of his
textbook \cite{Poincare1912} - two texts marking Poincar{\'e}'s desire to
show off his new probabilistic credo. But one has to realize that
Poincar{\'e} wanted to convince himself above all and this leads one to
ask the question Jean-Paul Pier ironically used as the title for his
paper \cite{Pier1996}: did Poincar{\'e} believe in the calculus of
probability or not? Without pretending to give a final answer, one may
however observe that Poincar{\'e} very honestly sought for a
demarcation of the zone where it seemed to him that using probability
theory did not create a major problem. Hence the attempt to tackle
some fundamental questions in order to go beyond the defects that
Bertrand had ironically illustrated with his famous paradoxes: When
is it legitimate to let randomness intervene? What definition can be
given of probability? Which mathematical techniques can be developed
in order to obtain useful tools for physics, in particular for the
kinetic theory of gases? Borel, as we shall see in the third section,
would later remember Poincar{\'e}'s position.

In his 1907 text, Poincar{\'e} accurately defined the legitimate way to
call upon the notion of randomness. He saw essentially three origins
for randomness: the ignorance of a very small cause that we cannot
know but which produces a very important effect (such as the so-called
Butterfly Effect), the complexity of the causes which prevents us from
giving any explanation other than a statistical one (as in the kinetic
theory of gases), the intervention of an unexpected cause that we have
neglected. This was not too far from the Laplacian conception, which
should not surprise us very much since Poincar{\'e}, born in 1854, was a
child of a century for which Laplace was a tutelary
figure. However, Poincar{\'e} knew well the accusations accumulated
against Laplace's theory and he proposed several ways to adjust it:
randomness, even if it is connected to our ignorance to a certain extent,
is not only that, and it is important to define the nature of the
connection between randomness and ignorance. The conventionalist
posture on which we have already commented naturally made things
easier, but Poincar{\'e} did not seek facility. As he wrote in 1899:

\begin{quote}`How shall we know that two possible cases are equally
  probable? Will it be by virtue of a convention? If we state an
  explicit convention at the beginning of each problem, everything
  will be fine; we will just apply the rules of arithmetic and
  algebra, and pursue the computation to the end, such that our result
  leaves no room for doubt. But, if we want to apply it anywhere, we
  will have to prove the legitimacy of our convention, and we will find
  ourselves in front of the difficulty we thought we had
  eluded.'\footnote{ `Comment saurons nous que deux cas possibles sont
    {\'e}galement probables? Sera-ce par une convention? Si nous pla\c
    cons au d{\'e}but de chaque probl{\`e}me une convention explicite, tout
    ira bien; nous n'aurons plus qu'{\`a} appliquer les r{\`e}gles de
    l'arithm{\'e}tique et de l'alg{\`e}bre et nous irons jusqu'au bout du
    calcul sans que notre r{\'e}sultat puisse laisser place au
    doute. Mais, si nous voulons en faire la moindre application, il
    faudra d{\'e}montrer que notre convention {\'e}tait l{\'e}gitime, et nous nous
    retrouverons en face de la difficult{\'e} que nous avions cru {\'e}luder.'
    (\cite{Poincare1899}, p.262)}\end{quote}

In a remarkable creative achievement, Poincar{\'e} forged a method
allowing the objectification of some probabilities. Using it, if one
considers for instance a casino roulette with alternate black and red
sectors, even without having the slightest idea of how it is put into
motion, one may show it is reasonable to suppose that after a large
number of turns, the probability that the ball stops in a red zone (or
a black zone) is equal to 1/2. There are thus situations where one can
go beyond the hazy Laplacian principle of (in)sufficient reason as a
necessary convention to fix the value of the probability. The profound
{\it method of arbitrary functions}, which is based on the hypothesis
that at the initial time the distribution of the position where the ball
stops is arbitrary and shows that this distribution reaches an
asymptotic equilibrium and tends towards the uniform distribution, was
certainly Poincar{\'e}'s most important invention in the domain of
probability and we shall see later the spectacular course it took.

To conclude this survey, let us mention that in 1906, when Poincar{\'e}
was completing his report for the Cour de Cassation in the Dreyfus
case, despite his place as the pre-eminent French authority in the
theory of probability, he was, as far as his scientific thought was
concerned, somewhat in the middle of the ford between a completely
deterministic description of the world and our modern conceptions
where randomness enters as a fundamental ingredient. Poincar{\'e} kept this
uncomfortable position until the end of his life. It is besides
noticeable that at the precise moment when Borel was - so strikingly -
taking over, Poincar{\'e} did not seem to have been
particularly interested in the enterprise of his young follower. He
showed similarly disinterest for the fortunate
experiments of the unfortunate Bachelier: he had written, it is true,
a benevolent report on his thesis \cite{Bachelier1900} and had
sometimes helped him to obtain grants, but the two men had no further scientific contact \cite{Courtault2002}. Even more
suprising, Poincar{\'e} seems to have thoroughly ignored the Russian
school's works (Chebyshev, Markov, Lyapounov\dots ) and this explains
why he became never conscious of some connections with his own
works. Poincar{\'e}'s probabilistic studies leave therefore a feeling of
incompleteness, partly due probably to his premature death at the age
of 58, but also to the singular situation of this last giant of
Newtonian-Laplacian science who remained on the threshold of upheavals that arrived after his death.

\section{Second part:  construction of a probabilistic approach}

In this second part, I would like to present some steps which have
marked the progressive entry of questions on probability in Poincar{\'e}'s
works. Even if it is not possible to see a perfect continuity in the
chain of his research, a kind of genealogy can be traced which allows
one to better understand how the mathematician gradually adopted a
probabilistic point of view in several situations. { I feel a need
  to mention that my aim in this section is to describe Poincar{\'e}'s
  contribution and not necessarily to comment its originality. A particular problem
  with our hero is that he practically never quotes his sources so
  that it is difficult to make a statement about what he had read or
  not.} I have tried as much as was possible to make each subsection
of this part independent of the others, which may sometimes result in
brief repetition.

\subsection{The recurrence theorem and its `probabilistic' extension}

In anticipation of the sixtieth birthday of King Oscar II of Sweden in
1889, a mathematical competition was organized by Mittag-Leffler. The
subject concerned the three-body problem: Was the system
Earth-Moon-Sun stable? Periodic? Organized so that it will always remain
in a finite zone of space? Many of these questions were fundamental
ones that had challenged Newtonian mechanics from the 18th
century. Poincar{\'e} submitted in 1888 an impressive memoir, immediately
selected for the award by a jury including Weierstrass, Mittag-Leffler and
Hermite. While correcting the proofs of the paper for {\it Acta
  Mathematica}, Phragmen located a mistake, leading Poincar{\'e} to make
numerous amendments before resubmitting a lengthy paper the following
year, published in volume 13 of {\it Acta Mathematica}.

This story, well known and well documented (see in particular
\cite{BarrowGreen1997}), is of interest for us only as far as one
difference between the version submitted for the prize and the
published version in 1890 is the appearance of the word {\it
  probability}, certainly for the first time in the French
mathematician's works. I shall closely follow Bru's investigation
\cite{Bru2013} on the way in which countable operations gradually
established themselves in the mathematics of randomness.

In the first part of the memoir submitted for the competition,
Poincar{\'e} studied the implications of the existence of integral
invariants on the behavior of dynamical systems. {A simple example of
  this situation is given by the case of an incompressible flow for
  which the shape of the set of molecules changes, but not its volume,
  which remains constant in time. However, Poincar{\'e} had in view a much
  more general situation. He showed in fact on p.46 {\it and seq.} of
  \cite{Poincare1890} that with a proper choice of coordinates in the
  phase space of position and velocity, the general situation of a
  mechanical system may be expressed by such an integral invariant.}
Poincar{\'e} then expounded a version of his recurrence theorem
in the following form: let $E$ be a bounded portion of the space,
composed of mobile points following the equations of mechanics, so
that the total volume remains invariant in time. Let us suppose,
writes Poincar{\'e}, that the mobile points remain always in $E$. Then, if
one considers $r_0$ a region of $E$, no matter however small it may
be, there will be some trajectories which will enter it an infinite
number of times.

Poincar{\'e}'s proof is a model of ingenuity and simplicity. One
discretizes time with a step of amplitude $\tau$. Let us call with
Poincar{\'e} $$r_1,r_2,\dots , r_n, \dots $$ the ``consequents'' of $r_0$,
that is to say, the successive positions of the different points of
the region $r_0$ at times $$\tau , 2\tau , \dots , n\tau , \dots $$ In
the same way, the ``antecedent'' of a region is the region of which it
is the immediate consequent. Each region $r_i$ has the same volume; as
they remain by hypothesis inside a bounded zone, some of them
necessarily intersect.  Let two such regions be $r_p$ and $r_q$ with
$p<q$, with intersection the region $s_1$ having nonzero volume: a
point starting in $s_1$ will be back in $s_1$ at time
$(q-p)\tau$. Going back in time, let us call $r_0^1$ the sub-region of
$r_0$ whose $p$-th consequent is $s_1$. A particle starting from
$r_0^1$ will again enter this region at time $(q-p)\tau$. We now start
again the process by replacing $r_0$ by $r_0^1$, and thus build a
decreasing sequence $(r_0^n)$ of sub-regions of $r_0$ such that each
point starting from $r_0^n$ returns $n$ times at least. Considering
a point in the intersection of the $r_0^n$ (whose non-emptiness is
taken for granted by Poincar{\'e}), a trajectory starting from such a
point will pass an infinite number of times through $r_0$.

In this form, the theorem is therefore completely deterministic. What
then compelled Poincar{\'e} to believe it necessary, in the new version
published in 1890, to rewrite his result in a probabilistic setting
(\cite{Poincare1890}, pp.71-72)? At first glance indeed, the
appearance of the word probability in these pages may seem surprising
in the chosen framework of the mechanics of Newton, Laplace and
Hamilton. In fact, as Bru remarks (\cite{Bru2013}), one must not be
misled by this use of probability by Poincar{\'e}, viewing it as a sudden
revelation of the presence of randomness having the ontic value we
spontaneously give it today. Poincar{\'e} himself wrote: {\it je me
  propose maintenant d'expliquer pourquoi [les trajectoires non
  r{\'e}currentes] peuvent {\^e}tre regard{\'e}es comme exceptionnelles}. What
Poincar{\'e} was thus looking for was a convenient way of expressing the
rarity, the thinness of a set. He was writing before the decisive
creation of the measure theoretic tools and particularly of Borel's
thesis which would, four years later, prove that a countable set has a
measure equal to zero. For a long while, astronomers in particular
employed the concept of probability with the meaning of practical
rarity and certainly one should not seek for a more sophisticated
explanation to justify the presence of the word coming from Poincar{\'e}'s
pen. It was a convenient way of speaking, whose aim was almost
entirely to hide the {\it obscure instinct} mentioned by Poincar{\'e} in
his 1899 paper (\cite{Poincare1899}, p.262), almost as an apology,
because we cannot do without it if we want to do scientific work.

Poincar{\'e} begins by expounding the following `definition': if one calls
$p_0$ the probability that the considered mobile point starts from a
region $r_0$ with volume $v_0$ and $p'_0$ the probability that it
starts from another region $r'_0$ with volume $v'_0$, then
$$
\frac{p_0}{p'_0}=\frac{v_0}{v'_0}.
$$ 
In particular, if $r_0$ is a region with volume $v$, used as a
reference, the probability that the mobile point starting from $r_0$
starts from a sub-region $\sigma_0$ with volume $w$ is given by
$\displaystyle\frac{w}{r v}$.  Equipped with this notion, the
mathematician wants to prove that the initial conditions in $r_0$ such
that the trajectory does not reenter $r_0$ more than $k$ times form a
set with probability zero, no matter how large the integer $k$.

Earlier in his paper, Poincar{\'e} had proved that if $r_0, \dots ,
r_{n-1}$ were $n$ regions with the same volume $v$ included in a
common region with volume $V$, and if $nv>kV$, then it was necessary
that there were at least $k+1$ regions whose intersection was
nonempty. Indeed, if one supposes that all the intersections taken
$k+1$ by $k+1$ were empty, one may write (in modern notation) that
$\displaystyle \sum_{i=0}^{n-1}\ind_{r_i}\le k$, hence $nv\le kV$ by
integrating over the volume $V$.

Let us still suppose valid the hypothesis of the previous theorem
which asserts that the mobile point remains in a bounded region, in a
portion of the space with volume $V$, and let us again take the
discrete step $\tau$ in time. Let us next choose $n$ sufficiently
large so that $n>\displaystyle \frac{kV}{v}$. One may then find, among
the $n$ successive consequents of a region $r_0$ with volume $v$,
$k+1$ consequents, denoted $$r_{\alpha_0}, r_{\alpha_1}, \dots ,
r_{\alpha_k}$$ with $\alpha_1<\alpha_2<\dots <\alpha_k$, having a
nonempty intersection denoted by $s_{\alpha_k}$. Let us now call $s_0$
the $\alpha_k$-th antecedent of $s_{\alpha_k}$ and $s_p$, the $p$-th
consequent of $s_0$. If a mobile point starts from $s_0$, it will
enter the regions $$s_0,s_{\alpha_k-\alpha_{k-1}},
s_{\alpha_k-\alpha_{k-2}}, \dots , s_{\alpha_k-\alpha_2},
s_{\alpha_k-\alpha_1}, s_{\alpha_k-\alpha_0}$$ which, by construction,
are all included in $r_0$ (as for each $0\le i\le k$, the
$\alpha_i$-th consequent of $s_{\alpha_k-\alpha_{i}}$ is in
$s_{\alpha_k}$ and therefore in $r_{\alpha_i}$) . One has therefore
shown that there are, in the considered region $r_0$, initial
conditions of trajectories which pass at least $k+1$ times through
$r_0$.

Let us eventually fix a region $r_0$ with volume $v$. Let us consider,
writes Poincar{\'e}, $\sigma_0$ the subset of $r_0$ such that the
trajectories issued from $\sigma_0$ do not pass through $r_0$ at least
$k+1$ times between time 0 and time $(n-1)\tau$; denote by $w$ the
volume of $\sigma_0$. The probability $p_k$ of the set of such
trajectories is therefore $\displaystyle \frac{w}{v}$.

By hypothesis a trajectory starting from $\sigma_0$ does not pass
$k+1$ times through $r_0$ , and hence not through $\sigma_0$. From the
previous result, one has necessarily that $nw<kV$, such that $$p_k<
\frac{kV}{nv}.$$ No matter how large $k$ may be, one may choose $n$
large so that this probability can be made as small as
wanted. Poincar{\'e}, tacitly using the continuity of probability along a
non-increasing sequence of events, concludes that the probability of
the trajectories issued from $r_0$ which do not pass through $r_0$
more than $k$ times between times 0 and $\infty$ is zero.

\subsection{Kinetic theory of gases}

As observed in the introduction, in 1892 Poincar{\'e} was not
favorably disposed towards the statistical description of
thermodynamics. His polemics with Tait, from which I quoted several
passages, was closely tied to the mechanist spirit in which Poincar{\'e}
had been educated. Statistical mechanics, and in particular the
kinetic theory of gases, could not therefore pretend to be more than
an ingenious construction with no explanatory value. An important
text revealing Poincar{\'e}'s thoughts on the subject was published
immediately afterwards, in 1893, in one of the first issues of the
{\it Revue de M{\'e}taphysique et de Morale} \cite{Poincare1893}.  With
great honesty, Poincar{\'e} mentioned the classical mechanical
conception of the universe due to Newton and Laplace, but also the
numerous problems it encounters when it tries to explain numerous
practical situations of irreversibility, as in the case of
molecular motion in thermodynamics. Poincar{\'e} mentioned that the
kinetic theory of gases proposed by the English is the attempt {\it la
  plus s{\'e}rieuse de conciliation entre le m{\'e}canisme et l'exp{\'e}rience}
(\cite{Poincare1893}, p.536). Nevertheless, he stated that numerous
difficulties still remained, in particular for reconciling the
recurrence of mechanical systems ({\it ``un th{\'e}or{\`e}me facile {\`a} {\'e}tablir''},
wrote the author, who may have adopted a humouristic posture) with
the experimental observation of convergence towards a stable
state. The manner in which the kinetic theory of gases pretends to
evacuate the problem, by invoking that what is called a stable
equilibrium is in fact a transitory state in which the system remains
an enormous time, did not seem to convince our hero. However, at least,
the tone adopted in \cite{Poincare1893} is obviously calmer than
that of the exchanges with Tait. Another point which can be observed is
that here, as in other works by Poincar{\'e} we shall comment on, Boltzmann was
nowhere to be found. This absence, difficult
to imagine as involuntary, remains unexplained, including for Von Plato
in \cite{Vonplato1991}, p.84.

In 1892, Lord Kelvin presented a note \cite{Kelvin1892} to the Royal
Society (of which he was then the president) with an unambiguous
title. The note presented an {\it ad hoc} example demonstrating, in a
supposedly {\it decisive} way, the failure of equipartition of kinetic
energy following Maxwell and Boltzmann's theory. The two physicists
had indeed deduced the equipartition of kinetic energy as a basic
principle of their theory: the average kinetic energies of several
independent parts of a system are in the same ratio as the ratio of
the number of degrees of freedom they have.  This result was
fundamental for the establishment of a relation between kinetic energy and temperature.

In his short paper, Kelvin imagined a mechanical system including
three points $A,B,C$, which are in motion in this order on a line
$KL$, such that $B$ remains almost motionless and only reacts to the
shocks produced by $A$ and $C$ on one side and the other, whereas the
mechanical situation on both sides is different because of a repulsive
force $F$ acting on $A$ and pushing it towards $B$ (in the zone $KH$
of the scheme) while $C$ can move freely.

\begin{figure}[ht!]
\begin{center}
 \includegraphics[height=10cm]{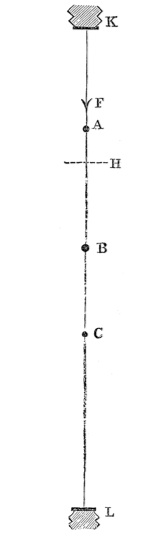}
 \caption{ Kelvin's construction in \cite{Kelvin1892}.}
  \end{center}
  \end{figure}

  The total energy of $C$ is balanced by the energy of $A$, but, as
  the latter includes a non-negative potential-energy term due to the
  repulsive force, Kelvin triumphantly concluded that the average
  kinetic energy of $A$ and $C$ cannot be equal, as they should have
  been following Maxwell's theory as the two points each have one
  degree of freedom. Kelvin commented:

  \begin{quote} `It is in truth only for an approximately ``perfect''
    gas, that is to say, an assemblage of molecules in which each
    molecule moves for comparatively long times in lines very
    approximately straight, and experiences changes of velocity and
    direction in comparatively very short times of collision, and it
    is only for the kinetic energy of the translatory motions of the
    molecules of the ``perfect gas'' that the temperature is equal to
    the average kinetic energy per molecule.' (\cite{Kelvin1892},
    p.399)\end{quote}

  Reading this note encouraged Poincar{\'e}, as he noted, to
  reflect on the kinetic theory of gases, to understand whether
  Kelvin's objection was well-founded, and to draw his own conclusions
  on the subject. At this precise moment, when he had been attacked by
  Tait, the title of the note by such an authority as Kelvin had
  impressed him and he may have thought he would find there a decisive
  argument confirming his own skepticism. And so in 1894 Poincar{\'e}
  published his first paper on the kinetic theory of gases
  \cite{Poincare1894}. Poincar{\'e} began by presenting a long general
  exposition of the fundamentals of Maxwell's theory. This survey
  seemed necessary in the first place because the kinetic theory of
  gases had been much less studied by French physicists than by
  English physicists\footnote{ `[\dots ] a {\'e}t{\'e} beaucoup moins cultiv{\'e}e
    par les physiciens fran\c cais que par les anglais.'
    (\cite{Poincare1894}, p.513)}. These fundamentals were in the first place the
  ergodic principle, called the {\it postulate of Maxwell} by
  Poincar{\'e}, which asserts that, whatever may be the initial situation
  of the system, it will always pass an infinity of times as close as
  desired to any position compatible with the integrals of the motion;
  from this postulate Maxwell drew a theorem whose main consequence
  was precisely the point contested by Kelvin: in a system for which
  the only integral is the conservation of the kinetic energy, if the
  system is made up of two independent parts, the long run mean values
  of the kinetic energy of these two parts are in the same ratio as
  their numbers of degrees of freedom.

  Poincar{\'e} began by observing, as previously in
  \cite{Poincare1893}, that the recurrence theorem of
  \cite{Poincare1890} contradicted Maxwell's postulate along the
  recurrent solutions. It was therefore necessary at least to add that
  the postulate was true except for certain initial
  conditions. \footnote{ `[\dots ] sauf pour certaines conditions
    initiales exceptionnelles' (\cite{Poincare1894}, p.518)} As von
  Plato comments \cite{Vonplato1991} (p.84), we have here the
  formulation usually given today for the ergodic principle, in order
  to take into account the possibility of exceptional initial
  conditions. Once again, although this idea was also present in
  Boltzmann's works, the Austrian scientist was nowhere mentioned.

  But it was above all the objection contained in Kelvin's paper that
  Poincar{\'e} desired to analyze in detail in order to check whether it
  contradicted Maxwell's results or not. From the situation of the
  system $A$, $C$, Poincar{\'e} built a representative geometric model: a
  point $M$ in a phase space with three dimensions whose first
  coordinate is the speed of $A$, the second the speed of $C$ and the
  third the abscissa of $A$. Using Kelvin's system conditions, he
  could define $S$, a solid of revolution, from which $M$ cannot exit
  in the course of time. Naturally, two small regions included in $S$
  with the same volume can be entered a different number of times with
  the same total sojourn time because the speed in these volumes could
  be different. Poincar{\'e} introduced the notion of the {\it density of
    the trajectory} in a small element in $S$ with volume $v$ as the
  quotient $\displaystyle \frac{t}{v}$, where $t$ is the total time
  spent by the trajectory in $v$ (\cite{Poincare1894}, p.519).  Using
  this representation, Poincar{\'e} could define the average value of the
  kinetic energy for $A$ as the moment of inertia of $S$ with respect
  to the plane $yz$, for $C$ as the moment of inertia with respect to
  $xz$, the `masses' in $S$ being distributed by the previously
  defined density. The solid $S$ being one of revolution, these
  moments of inertia are equal: the fine analysis made by Poincar{\'e}
  therefore shows that one can recover the equipartition result by
  taking the average of the kinetic energies not uniformly over time
  but taking into account the phases of the motion and their duration.

  Poincar{\'e} concluded his paper with a comment that may seem
  paradoxical in light of the result he had just obtained. While
  disputing the decisive character of Kelvin's arguments, Poincar{\'e}
  insisted that he nevertheless shared his colleague's skepticism. To
  give weight to his comment, he slightly transformed Kelvin's example
  in order to produce an {\it ad hoc} situation for which there {\it
    really is} a problem. In fact, some lines earlier, Poincar{\'e} had
  emphasized what was for him the fundamental point:

  \begin{quote} `I believe Maxwell's theorem really is a
    necessary consequence of his postulate, as soon as one admits the
    existence of a mean state; but the postulate itself must admit
    many exceptions.'\footnote{ `Je crois que le th{\'e}or{\`e}me de Maxwell
      est bien une cons{\'e}quence n{\'e}cessaire de son postulat, du moment
      qu'on admet l'existence d'un {\'e}tat moyen; mais le postulat
      lui-m{\^e}me doit comporter de nombreuses exceptions.'
      (\cite{Poincare1894}, p.521)}\end{quote}

  For Poincar{\'e}, it was therefore the definition of the average
  states that posed a problem, and it was the search for a
  satisfactory definition that required attention from those those who wished to
  consolidate the bases of statistical mechanics. We shall
  see later that it was indeed in this direction that Poincar{\'e}, and
  later Borel, focussed their efforts.

\subsection{Limit theorems}

The textbook \cite{Poincare1896} published in 1896 constitutes the
first of Poincar{\'e}'s works dealing explicitly with the theory of
probabilities. It was, as already mentioned, the result of
lectures offered by Poincar{\'e} during the academic year 1893-94 at the
Sorbonne, redacted by a former student of the {\'E}cole Normale Sup{\'e}rieure,
who became an actuary and who probably wished to learn with the
master; it was published by Georges Carr{\'e}. This first edition does not
have a preface, and presents itself as a succession of 22 lectures,
more or less connected to each other, probably reflecting Poincar{\'e}'s
actual lectures. It is 274 pages long, compared to the 341 pages of
the second edition of 1912 \cite{Poincare1912}, which provides an idea of the
number of complements. In its initial form,
Poincar{\'e}'s book appears as a successor to Bertrand's treatise
\cite{Bertrand1888}, the framework of which it follows. Nevertheless,
these textbooks really have a different feel, and we must leave it at
that for the moment. However, as Poincar{\'e}'s textbook has been
discussed by several commetators, especially in
\cite{Sheynin1991} and \cite{Cartier2010}, I shall restrict myself to
just a few remarks. Let us note that the authors just mentioned
focused on the 1912 edition, which naturally benefited from
reflections of Poincar{\'e} after that very important period of the
mid-1890s, when he was beginning to investigate the mathematics of
randomness, so that this may not accurately reflect the
mathematician's state of mind in 1896. I choose here, in
contrast, to focus on the original 1896 version.

An important part of Poincar{\'e}'s textbook is devoted to the use of
probability theory as a model of measurement error in the experimental
sciences. In a commentary on his own works (\cite{Poincare1921},
p.121), Poincar{\'e} wrote:

\begin{quote}`The Mathematical Physics Chair has for its official
  title: Calculus of Probability and Mathematical Physics. This
  connection can be justified by the applications of this calculus in
  all physics experiments, or by those in kinetic gas theory. In any
  case, for one semester I dealt with probability and my lectures were
  published. The theory of errors was naturally my main goal. I had to
  erect explicit reservations concerning the generality of the `law of
  errors', but I sought to justify the law via new considerations in
  the cases where it remains legitimate.'\footnote{ `La Chaire de
    Physique Math{\'e}matique a pour titre officiel\,: Calcul des
    Probabilit{\'e}s et Physique Math{\'e}matique. Ce rattachement peut se
    justifier par les applications que peut avoir ce calcul dans
    toutes les exp{\'e}riences de Physique; ou par celles qu'il a trouv{\'e}es
    dans la th{\'e}orie cin{\'e}tique des gaz. Quoi qu'il en soit, je me suis
    occup{\'e} des probabilit{\'e}s pendant un semestre et mes le\c cons ont
    {\'e}t{\'e} publi{\'e}es. La th{\'e}orie des erreurs {\'e}tait naturellement mon
    principal but. J'ai d{\^u} faire d'expresses r{\'e}serves sur la
    g{\'e}n{\'e}ralit{\'e} de la ``loi des erreurs"; mais j'ai cherch{\'e} {\`a} la
    justifier, dans les cas o{\`u} elle reste l{\'e}gitime, par des
    consid{\'e}rations nouvelles.'}\end{quote}

In \cite{Poincare1896}, the analysis of the law of errors begins on
page 147 and occupies part of the following chapters. Poincar{\'e}
commented on the manner in which the Gaussian character of
the error had been obtained:

\begin{quote}`[This distribution] cannot be obtained by rigorous
  deductions. Several of its putative proofs are awful, including,
  among others, the one based on the statement that gap probability is
  proportional to the gaps. Nonetheless, everyone believes it, as
  M. Lippmann told me one day, because experimenters imagine it to be
  a mathematical theorem, while mathematicians imagine it to be an
  experimental fact.'\footnote{`[Cette loi] ne s'obtient pas par des
    d{\'e}ductions rigoureuses; plus d'une d{\'e}monstration qu'on a voulu en
    donner est grossi{\`e}re, entre autres celle qui s'appuie sur
    l'affirmation que la probabilit{\'e} des {\'e}carts est proportionnelle
    aux {\'e}carts. Tout le monde y croit cependant, me disait un jour
    M.~Lippmann, car les exp{\'e}rimentateurs s'imaginent que c'est un
    th{\'e}or{\`e}me de math{\'e}matiques, et les math{\'e}maticiens que c'est un fait
    exp{\'e}rimental.'(\cite{Poincare1896}, p.149)}\end{quote}

{Much of Poincar{\'e}'s treatment of the law of errors, in particular the
  convergence to a Gaussian distribution of a Bayesian approach of the
  measurement process, was in fact already found by Laplace, at
  least in a preliminary form. Once again, it is not really possible to
  make precise what Poincar{\'e} knew. For details on what is
  today known as the Bernstein-von Mises theorem see \cite{Vaart1998}
  (p.140 {\it et seq.} and its references).}

Let us consider observations of a phenomenon denoted by $x_1,x_2,
\dots , x_n$. The true measure of the phenomenon under study being
$z$, the {\it a priori} probability that each of these $n$
observations belong to the interval $[x_i,x_i+dx_i]$ is taken under
the form
$$\varphi (x_1,z)\varphi (x_2,z)\dots \varphi (x_n,z)dx_1dx_2\dots dx_n.$$

Let finally $\psi (z)dz$ be the {\it a priori} probability so that the
true value belongs to the interval $[z,z+dz[$.

Supposing that $\psi$ is constant and that $\varphi (x_i,z)$ can be
written under the form $\varphi (z-x_i)$, Gauss obtained the
Gaussian distribution by looking for the $\varphi$ such that the most
probable value was the empirical mean $$\overline{x}=\frac {x_1+\dots
  +x_n}{n}.$$ Poincar{\'e} recalled (\cite{Poincare1896}, p.152)
Bertrand's objections to Gauss' result; Bertrand had in particular
disputed the requirement that the mean be the most probable value
while the natural condition would have been to require it to be the
{\it probable value} (which is to say the expectation).

Poincar{\'e} thus considered the possibility of suppressing some of
Gauss' conditions. Keeping firstly the hypothesis that the empirical
mean be the most probable value (\cite{Poincare1896}, p.155 - see also
the details in \cite{Sheynin1991}, p.149 {\it et seq.}), he obtained
for the form of the error function
$$\varphi (x_1,z)=\theta (x_1)e^{A(z)x_1+B(z)},$$
where $\theta$ and $A$ are two arbitrary functions, $B$ being such that the following differential equation is satisfied: $A'(z) z+B'(z)=0$.   

Considering next Bertrand's objection, Poincar{\'e} looked next at the
problem that arises when one replaces the requirement of the {\it most
  probable value} with that of the {\it probable value}. Here
(\cite{Poincare1896}, p.158) he gives a theorem he would use
subsequently a number of times: if $\varphi_1$ and $\varphi_2$ are two
continuous functions, the quotient
$$\frac{\int \varphi_1(z)\Phi^p(z)dz }{ \int \varphi_2(z)\Phi^p(z)dz}$$ tends, when  $p\rightarrow +\infty$, towards 
$$\frac{\varphi_1(z_0)}{\varphi_2(z_0)},$$ where $z_0$ is a point in
which $\Phi$ attains its unique maximum. Following his habit, which
drove Mittag-Leffler to dispair, Poincar{\'e}'s writing was somewhat
laconic; he did not specify any precise hypothesis, or supply a real proof, presenting the result only as an extrapolation of the discrete case.

In any case, considering next $$\Phi (x_1,\dots ,x_n; z)= \varphi
(x_1,z) \varphi (x_2,z) \dots \varphi (x_n,z),$$ Poincar{\'e} made the
hypothesis that $p$ observations resulted in the value $x_1$, $p$
resulted in $x_2$, \dots , $p$ resulted in $x_n$, where $p$ is a fixed
and very large integer (\cite{Poincare1896}, p.157).

The condition requiring the mean to be equal to the expectation can therefore be written as 

$$
\frac{\int_{+\infty}^{+\infty}z\psi (z)\Phi^p (x_1,\dots ,x_n; z)dz}{ \int_{+\infty}^{+\infty}\psi (z)\Phi^p (x_1,\dots ,x_n; z)dz}=\frac{x_1+\dots +x_n}{n}. 
$$
Applying the previous theorem, under the hypothesis that $\Phi$ has a
unique maximum in $z_0$, one has $z_0$ as the limit of the left-hand
side, which must therefore be equal to the arithmetic mean $\overline
x$. One is thus brought back to the previous question under the
hypothesis that $\Phi$ should be maximal at $\overline x$. Under the
hypothesis that $\varphi$ depends only on the discrepancies $z-x_i$,
Poincar{\'e} again obtained the Gaussian distribution. It is remarkable
that the form of the {\it a priori} probability of the phenomenon
$\psi$ is not present in the result. This lack of dependence on the
initial hypothesis might perhaps have been the inspiration for his
method of arbitrary functions, described later.

Poincar{\'e} examined next the general problem by suppressing the
constraint that $\varphi$ depends only on the discrepancies, and
obtaind the following form for $\varphi$
$$\varphi (x_1,z)=\theta (x_1)e^{-\int \psi (z)(z-x_1)dz}$$ where $\int \psi (z)(z-x_1)dz$ is the primitive of $\psi (z)(z-x_1)$ equal to 0 in $x_1$.

He argued (\cite{Poincare1896}, p.165) that the only reasonable
hypothesis was to take $\psi=1$ as there was no reason to believe that
the function $\varphi$, which depends on the observer's skillfulness,
would depend on $\psi$, the {\it a priori} probability for the value
of the measured quantity. For $\theta$, on the contrary, there was no
good reason to suppose it constant (in which case the Gaussian
distribution would be again obtained). Poincar{\'e} took the example of
the meridian observations in astronomy where a {\it decimal error} had
been detected in practice: the observers show a kind of predilection
for certain decimals in the approximations.

Poincar{\'e} gave a somewhat intricate justification for focusing on the
mean because it satisfies a practical aspect: as the errors are small,
to estimate $f(z)$ by the mean of the $f(x_i)$ was the same as
estimating $z$ by the mean of the $x_i$, as immediately seen by
replacing $f(x)$ by its finite Taylor expansion in
$z$, $$f(z)+(x-z)f'(z).$$

In any case, the major justification was given in the following
chapter (Quatorzi{\`e}me le\c con, \cite{Poincare1896}, p.167), where the
consistency of the estimator $\overline x$ was studied using an
arbitrary law for the error, based on the law of large numbers. After
having recalled the computation of the moments for the Gaussian
distribution, Poincar{\'e} implemented the method of moments in the
following way. Suppose that $y$, with distribution $\varphi$, admits
the same moments as the Gaussian distribution. One then computes the
probable value of $e^{-n(y_0-y)^2}$, where $n$ is a given integer.
Decomposing $$e^{-n(y_0-y)^2}=\sum_{p=0}^\infty A_py^{2p},$$ one
obtains $$\int_{-\infty}^{+\infty}\sqrt{h/\pi
}e^{-hy^2}e^{-n(y_0-y)^2}dy=\sum_{p=0}^\infty A_p\, \mathbb
E(y^{2p}),$$ letting $\mathbb E(y^{2p})$ denote the expectation of
$y^{2p}$ (the odd moments are naturally equal to zero), and $h$ a
positive constant. The same decomposition is valid by hypothesis if
$\varphi$ replaces the Gaussian distribution in the integral. One has
therefore
$$
\frac{\int_{-\infty}^{+\infty}\sqrt{h/\pi }e^{-hy^2}e^{-n(y_0-y)^2}dy}{\int_{-\infty}^{+\infty}\varphi (y)e^{-n(y_0-y)^2}dy}=1,
$$ 
and, using his theorem on the limits again, Poincar{\'e} could obtain, letting $n$ tend to infinity, $$\sqrt{h/\pi }e^{-hy_0^2}= \varphi (y_0).$$

Now, let us consider again that $n$ measures $$x_1,\dots , x_n$$ of a
quantity $z$ are effectuated, and let us denote by $y_i=z-x_i$ the
individual error of the $i$-th measure. Let us suppose that the
distribution of an individual error is arbitrary.

Poincar{\'e} began by justifying the fact of considering the mean 
$$
\frac{y_1+\dots +y_n}{n}
$$  
of the $n$ individual errors as error. Indeed, he explained that the mean becomes more and more probable as the probable value of its square is  
$$\frac{1}{n}\mathbb E ({y_1^2}),$$  
and so, when $n$ becomes large, the probable value of 
$$\left(\frac{y_1+\dots +y_n}{n}\right)^2$$ 
tends towards 0 in the sense that the expectation 
$$\mathbb E[(z-\overline x)^2]$$ 
tends towards 0 (this is the $L^2$ version of the law of large
numbers). As Sheynin observes (\cite{Sheynin1991}, p. 151), Poincar{\'e}
made a mistake when he attributed to Gauss this observation as the
latter had never been interested in the asymptotic study of the
error. Poincar{\'e} then used his method of moments in order to
prove that the distribution of the mean is Gaussian when the
individual errors are centered and do not have a significant effect on
it.

In the 1912 edition \cite{Poincare1912}, Poincar{\'e} significantly added
a section (\cite{Poincare1912}, n$^\circ$144 pp. 206-208) devoted to a
proof of the central limit theorem and obtained a {\it justification
  {\rm a posteriori} de la loi de Gauss fond{\'e}e sur le th{\'e}or{\`e}me de
  Bernoulli}. Poincar{\'e} introduced the characteristic function
as $$f(\alpha )=\sum_x p_x e^{\alpha x}$$ in the {\it finite} discrete
case, and $$f(\alpha )=\int \varphi (x)e^{\alpha x}dx$$ in the case of
a continuous density. In his mind, $\alpha$ was a real or a complex
number and neither the bounds of the sum or the integral, nor the
issue of convergence were mentioned. { Poincar{\'e}'s considerations on
  the Fourier and Laplace transforms can in fact already be largely
  found in Laplace, though the denomination {\it characteristic
    function} is Poincar{\'e}'s. The Fourier inversion formula was used by
  Poincar{\'e} since his lectures on the analytical theory of heat
  \cite{Poincare1895} (see in particular the Chapter 6, p.97). As
  usual, there is no mention of Laplace in Poincar{\'e}'s text. Plausibly,
  Laplace introduced characteristic functions for probability
  distributions after having studied Fourier's treatise (see
  \cite{Bru2012}). In the absence of proof that Poincar{\'e} knew
  Laplace's method, one may speculate that he had the same kind of
  illumination as his predecessor.}

Thanks to Fourier's inversion formula, Poincar{\'e} stated that the
characteristic function determined the distribution. He could thus
obtain simply that a sum of independent Gaussian variables followed a
Gaussian distribution and, by means of a heuristic and again quite
laconic proof, that the error resulting from a large number of very
small and independent partial errors\footnote{ `[\dots ] r{\'e}sultante
  d'un tr{\`e}s grand nombre d'erreurs partielles tr{\`e}s petites et
  ind{\'e}pendantes.' (\cite{Poincare1912}, p.208)} was Gaussian. It seems
difficult to award the status of a proof of the central limit theorem
to these few lines, a proof published some ten years later in works by Lindeberg and L{\'e}vy
(\cite{Lindeberg1922}, \cite{Levy1922}). Besides, in this intriguing
but rather hasty complement, Poincar{\'e} showed his complete ignorance of
the Russian research on limit theorems (\v Ceby\v cev, Markov et
Lyapunov) which produced some well established versions of the
theorem.

In the sixteenth chapter of \cite{Poincare1896}, (n$^\circ$147 of the
second edition \cite{Poincare1912} p.211), Poincar{\'e} the physicist
still had reservations about what would be an indiscriminate use of
the theories he had just described, which depended so heavily on a
mathematical idealization (the absence of systematic errors, overly
smooth hypotheses\dots ). He wrote, not without irony: `I've argued the
best I could up to now in favor of Gauss' law.'\footnote{ `J'ai
  plaid{\'e} de mon mieux jusqu'ici en faveur de la loi de Gauss.'}  He
then focused on the study of exceptional cases and completed his
textbook wit a detailed examination of the method of least squares; on
these subjects, I refer the interested reader to the already quoted
paper by Sheynin (\cite{Sheynin1991}).

\subsection{The great invention: the method of arbitrary functions}

Although Poincar{\'e} is rightly considered to be the father of
conventionalism in scientific philosophy, it would be simplistic to
think that this position covers all Poincar{\'e}'s philosophy of
research. Admittedly, from the very beginning, the latter had repeated
that any use of probability must be based on the choice of a
convention one had to justify. Thus, if one throws a die, one is
generally led to take as a convention the attribution of a probability
of 1/6 for each face to appear. However all the arguments used to
justify this convention do not have the same value, and choosing well
among them is also part of a sound scientific process.  We cannot indeed
be satisfied with common sense: Bertrand amused himself, when he
constructed his famous paradoxes about the choice of a chord in a
circle, in showing that the result depended so closely on the chosen
convention that it lost meaning, and that the calculus of probability in
such situations was reduced to a more or less ingenious arithmetic. The
risk was to condemn the calculus of probability altogether as a vain
science and conclude that our {\it instinct obscur} had deceived us
(\cite{Poincare1899}, p.262).

And yet, wrote Poincar{\'e}, without this obscure instinct science would
be impossible. How can one reconcile the irreconcilable?

Until then, common practice was to cite Laplace and use the
principle of insufficient reason as a supporting argument. A dubious
argument in fact, as in practice it amounts to assigning a value to the
probability only by supposing that the different possible cases are
equally probable, since we do not have any reason to assert the
contrary. How could a scientist such as Poincar{\'e}, who was looking for
a reasonably sound basis for using the mathematics of randomness, be
satisfied with such a vicious circle?

Let us observe in passing that he was far from the first to deal
with such a question. And besides, we have already seen that after
Laplace's death, several weaknesses of his approach had been
underlined: the vicious circle of the definition of probability by
possibility, the absence of an answer to the general question of the
nature of the probabilities of causes when applying Bayes's principle,
let alone the confusions in ill-considered applications, in particular
the judicial ones that we have already mentioned\dots. A substitute was
sought for Laplace's theory. This problem of defining the {\it
  natural} value of probabilities had in particular obsessed German
psychologists and physiologists throughout the second half of the 19th
century (\cite{Kamlah1987}). Von Kries in particular succeeded, a good
ten years before Poincar{\'e}, in constructing the foundations of a method
allowing one to justify the attribution of equal probabilities to the
different outcomes of a random experiment repeated a large number of
times (\cite{Kamlah1983}). Poincar{\'e}, without question, completely
ignored these works, all the more because they did not belong {\it
  stricto sensu} to the sphere of mathematics.

The question thus for Poincar{\'e} was to show that in some important
cases, one may consider that the equiprobability of the issues in a
random experiment as the result not only of common sense but also of
mathematical reasoning, and thereby to avoid the criticism of Laplace's
principle.

The idea developed by Poincar{\'e}, as earlier by von Kries, was that the
repetition of the experiment a large number of times ends in a kind of
asymptotic equilibrium, in a compensation, so that the hypothesis of
equiprobability becomes reasonable even if one absolutely ignores the
initial situation.

{ As early as the 1780s, Laplace observed that in many
  cases the initial distribution asymptotically vanishes when one
  repeatedly applies Bayes' method. For complements, the interested
  reader may consult Stigler's translation and commentary of Laplace
  (\cite{Stigler1986}), and also \cite{Bru2003}, p.144.  At the end of
  his paper \cite{Poincare1891}, Poincar{\'e} briefly mentioned that
  the convention he had adopted to define the probability concept in
  the study of the stability of the three-body problem was by no means
  necessary: the result (a probability equal to zero for non-recurrent
  trajectories) remains true whatever the convention. Poincar{\'e}'s aim in subsequent papers was to show that in
  more general situations, the knowledge of the probability
  distribution at the origin of time is not needed, as this
  distribution is not present in the final result. Two examples serve
  repeatedly as an illustration in the various texts (\cite{Poincare1899},
  \cite{Poincare1907} especially) in which he discussed his method: the
  uniform distribution of the so-called small planets on the Zodiac
  and the probabilities for the red and black slots of a roulette
  wheel.  In \cite{Poincare1899}, p. 266, Poincar{\'e} observed the
  proximity of these two situations. The systematic exposition of {\it
    arbitrary functions} as a fundamental method of the theory of
  probability seems due to Borel (\cite{Borel1909}, p.114) but the
  expression {\it m{\'e}thode des fonctions arbitraires} was later
  generally used in the context of Markov chains (see in particular
  \cite{Frechet1938}).}

Let us first follow Poincar{\'e}'s comments on the second, and simpler
case of roulette (\cite{Poincare1899}, p.267). The ball, thrown with
force, stops after having turned many times around the face of a
roulette wheel regularly divided into black and red sectors. How can
we estimate the probability that it stops in a red sector?

Poincar{\'e}'s idea is that, when the ball travels for a large number of
turns before stopping, any infinitesimal variation in the initial
impulsion can produce a change in the color of the sector where the
ball stops. Therefore, the situation becomes the same as considering
that the face of the game is divided into a large number of red and
black sectors. I make, said Poincar{\'e}, the convention that the
probability for this angle to be fall between $\theta$ and $\theta
+d\theta$ equals $\varphi(\theta )d\theta$, where $\varphi$ is a
function about which I do not know anything (as it depends on the way
the ball moved at the origin of time, an {\it arbitrary}
function). Poincar{\'e} nevertheless asserts, without any real
justification, that we are {\it naturally} led to suppose $\varphi$ is
continuous. The probability that the ball stops in a red sector is the
integral of $\varphi$ estimated on the red sectors.

Let us denote by $\varepsilon$ the length of a sector on the
circumference, and let us consider a double interval with length
$2\varepsilon$ containing a red and a black sector. Let then $M$ and
$m$ be respectively the maximum and the minimum of $\varphi$ on the
considered double interval. As we can suppose that $\varepsilon$ is
very small, the difference $M-m$ is very small. And as the difference
between the integral on the red sectors and the integral on the black
sectors is dominated by
$$\sum_{k=1}^{\pi / \varepsilon} (M_k-m_k)\varepsilon ,$$ (where $M_k$ and $m_k$ are respectively the maximum and the minimum on each double interval $k$ of the subdivision of the face with length $2\varepsilon$), 
this difference is small and it is thus reasonable to suppose that
both integrals, whose sum equals 1, are equal to 1/2.

Once again, Poincar{\'e}'s writing is somewhat sloppy. He emphasized the
importance of the fact that $\varepsilon$ was small with respect to
the total swept angle,\footnote{ `[\dots ] par rapport {\`a} l'angle total
  parcouru.'} but without giving much detail on how to interpret this
fact. His brevity probably came from the parallel with the other
example we shall now present - but which he studied earlier in his text.

The expression {\it small planets} designates the astero{\"\i}d belt
present between Mars and Jupiter which had been gradually explored
until the end of the 19th century. The first appearance of questions
of a statistical type about these planets seems to go back to the
twelfth chapter of the 1896 textbook (\cite{Poincare1896}, p. 142),
where Poincar{\'e} asked how one can estimate the probable value of their
number $N$. For that purpose, he implemented a Bayesian method using
the {\it a priori} probability for an existing small planet to have
been observed, this probability being supposed to have a density
$f$. It allowed him to carry out the computation of the {\it a
  posteriori} expectation of $N$.

In \cite{Poincare1899}, Poincar{\'e} was interested in a remarkable
phenomenon: the almost uniform distribution of the small planets in
the different directions of the Zodiac. Poincar{\'e} looked for arguments
justifying this fact (\cite{Poincare1899}, p. 265 {\it et seq.}). We
know, said Poincar{\'e}, that the small planets follow Kepler's laws, but
on the contrary we absolutely ignore their initial
distribution.

Let then $b$ be the longitude of a small planet at the initial time,
and $a$ its mean motion. At time $t$, its longitude is therefore
$at+b$.  As already noted, one knows nothing of the initial
distribution and we suppose it is given by an arbitrary function
$\varphi (a,b)$, once more assumed regular in some way: Poincar{\'e} wrote
that it was
{\it continuous} but in the sequel used it as a function of class
$C^\infty$.

The mean value of $\sin (at+b)$ is given by $$\int\int \varphi
(a,b)\sin (at+b)dadb.$$ When $t$ becomes large, this integral becomes
close to 0. Poincar{\'e} used successive integrations by part,
with the derivatives of $\varphi$, whereas he could have used only
continuity and Riemann-Lebesgue's lemma, but, as we have already seen,
Poincar{\'e} did not regard the refinement of his hypotheses as a major
concern.  {\it A fortiori}, for every nonzero integer $n$, the
integrals $$\int\int \varphi (a,b)\sin n(at+b)\,dadb,\,\,\,
\textrm{and} \,\,\,\int\int \varphi (a,b)\cos n(at+b)\,dadb$$ are also
very small for a large fixed $t$. Therefore, if one denotes by $\psi$
the probability density of the longitude at time $t$, one has for
every $n\ge 1$, $$\int_{[0,2\pi [}\psi (u)\sin nu\,
du,\,\,\,\textrm{and}\,\,\,\int_{[0,2\pi [}\psi (u)\cos nu\, du$$ very
close to 0. The Fourier expansion of $\psi$ leads to the conclusion
that $\psi$ is almost constant, that is to say, that the longitude of
a small planet is roughly uniformly distributed on the Zodiac.

\subsection{Cards shuffling}

If the example of the small planets illustrates sensitivity to initial
conditions, that of kinetic gas theory is connected with the
complexity of causes. The number of molecules is so large, and they
collide in so many ways, that it is impossible to consider the system
they form as simply describable by classical mechanics. In 1902 the
first textbook was published expounding the basic principles of
statistical mechanics, written by Gibbs (\cite{Gibbs1902}).  It
developed two main applications for the new theory: in addition to the
kinetic theory of gases, it introduced the situation of mixing two
liquids (a drop of ink put in a glass of water) in order to present
the evolution of a system towards equilibrium. Hadamard, in 1906, had
written a review of Gibb's book for the {\it Bulletin des sciences
  math{\'e}matiques} (\cite{Hadamard1906}). In order to illustrate this
mixing situation, he invented the ingenious metaphor of the shuffling
of a pack of cards by a gambler evolving towards an equal distribution
of the possible permutations of the cards. Hadamard however did not
propose any mathematical treatment of the question and it was
Poincar{\'e}, in the paper he published in 1907 in Borel's journal
(\cite{Poincare1907} republished later in \cite{Poincare1908} and
\cite{Poincare1912}), which first analyzed the problem. He restricted
himself in fact to the simplest case, that of two cards. Let us
suppose, said Poincar{\'e}, that one has a probability $p$ that after one
permutation, the cards are still in the same order as before the
permutation, and $q=1-p$ that their order is reversed. Let us consider
there are $n$ successive permutations and that the gambler who
shuffles the cards earns a payoff $S$ equal to 1 franc if the order
after these $n$ permutations is unchanged, and -1 franc if it is
reversed. A direct computation of the expectation shows that $$\mathbb
E(S) = (p-q)^n,$$ as, in a modern formulation, $S$ can be written
as $$\prod_{i=1}^n X_i$$ with the $X_i = \pm 1$ independent with
distribution $(p,q)$ representing the fact that the $i$-th permutation
has changed or not the order of the cards. Hence, except in the
trivial cases $p = 0$ or 1, $\mathbb E(S)\rightarrow 0$ when $n$ tends
towards infinity, which amounts to saying that the two states +1 and
-1, and therefore the two possible orders, tend to become
equiprobable. It is interesting that Poincar{\'e} was inspired, for
the recursive computation of the expectation without first looking for
the distribution of $S$, by several computations of expectations he
found in Chapter III of Bertrand's book \cite{Bertrand1888}.

As mentioned by Poincar{\'e}, the tendency to uniformity remained true
whatever the number of cards but the {\it d{\'e}monstration serait
  compliqu{\'e}e}. { One may suppose that in 1907 Poincar{\'e} already had the
  idea of the proof for the general case of $n$ cards but he wrote it
  only for the second edition of his textbook of probability in 1912,
  in the first section of a chapter added to the book, entitled {\it
    Questions diverses}. }  Curiously, Poincar{\'e}'s method of proof,
contrary to what he had done in the case of two cards, was not
inspired by probabilistic reasoning but was connected to the theory of
groups. { Though the Perron-Frobenius theory was already available
  (but probably unknown to Poincar{\'e}), Poincar{\'e} referred to older works
  by Frobenius published in the {\it Sitzungsberichte} of Berlin
  Academy between 1896 and 1901, and by Elie Cartan \cite{Cartan1898}
  that he had himself extended in his paper \cite{Poincare1903}} (consult
\cite{Seneta1981} and \cite{Cartier2010} for details). We shall see in
the next part that this non-probabilistic aspect did not escape Borel,
who proposed an alternative approach. As for card shuffling, it
enjoyed a spectacular renewed interest in the 1920s.

\section{Third part: an uneven heritage}

We now tackle the heritage of Poincar{\'e}'s ideas about randomness and
probability. This is an intricate question. Indeed, Poincar{\'e} cannot be
considered as a full probabilist in the spirit of mathematicians of
later generations like Paul L{\'e}vy and Andrei Nikolaevich Kolmogorov. As
we have already mentioned, these studies on probability constitute
only a very small island in the ocean of the mathematician's
production. Moreover, it is rather difficult to locate a very precise
result, a theorem, concerning the theory of probability which can be
specifically credited to Poincar{\'e}. His primary goal was to refine
already existing results or to explore new aspects and new questions
while not feeling compelled to give them a complete structure.  We
should again repeat here that in this domain more than any other, it
was above all Poincar{\'e} whom Poincar{\'e} sought to convince, and therefore
his works dealing with probability, including his philosophical texts,
often take on a rambling tone, written following his train of thought,
often slightly verbose, illustrating Picard's opinion (as reported by
K.~Popoff): he ignored the adage {\it pauca sed matura}\footnote{ `Il
  ne connaissait pas l'adage {\it pauca sed matura}.'
  (\cite{Popoff1993}, p. 89)}. As Bru notes (\cite{Bru2003}, p. 155),
everyone at that time had read Poincar{\'e}. But one has the impression
that few understood what he wrote about probability.

\subsection{Borelian path}

{\'E}mile Borel was unquestionably the main exception. Not only did he
read and understand Poincar{\'e}, but he was poised to make the subject his own in a
spectacular way, so much so that he may be regarded as the first
French probabilist of the twentieth century. We shall examine how this
passing of the baton took place between the master and his young
disciple.

It must firstly be said that this probabilistic turn of {\'E}mile Borel
was one of the most singular changes one could observe in a
mathematician around 1900. After initiating a profound
transformation of methods in the theory of functions, Borel became
a star of mathematical analysis in France. Nothing seemed to
predispose him to take the plunge and to devote important efforts to
study, refine and popularize the calculus of probability whose dubious
reputation in the mathematical community - on which we commented above
- might have led to rebuke. The context of Borel's turn since
1905, the date of publication of his first work in this domain, has
been studied in detail several times, for instance in the papers
\cite{DurandMazliak2011} and \cite{MazliakSage2012}. The main
difference which can be found between the discovery of probability by
Poincar{\'e} and by Borel is that, for the latter, it arose from
reflections within the mathematical field and more specifically by
considerations on the status of mathematical objects - in particular
about real numbers. In Borel, during the years just preceding 1900,
we note indeed a greater and greater distance from Cantorian
romanticism and its absolutist attitude, as emphasized by Anne-Marie
D{\'e}caillot in her beautiful book on Cantor and France
(\cite{Decaillot2008}, p. 159). Borel gradually replaced this
idealistic vision, which no longer satisfied him, by a realism colored
by a healthy dose of pragmatism: the probabilistic approach appeared
then to Borel as an adequate means by which to confront various forms of
reality: mathematical reality first and then physical reality and
practical reality\dots.

The best synthesis summing up Borel's spirit about the quantification
of randomness can be found in Cavaill{\`e}s' text \cite{Cavailles1940}
published in the {\it Revue de M{\'e}taphysique et de Morale}; it should
be seen, at least in part, as a commentary on Borel's essay on
the interpretation of probabilities \cite{Borel1939} that completed
the great enterprise of the {\it Trait{\'e} du Calcul des Probabilit{\'e}s et
  de ses Applications} begun in 1922. As Cavaill{\`e}s lyrically put it
(\cite{Cavailles1940}, p.154), probabilities appear to be the only
privileged access to the path of the future in a world which is no
longer equipped with the sharp edges of certainty, but presents
itself instead as a hazy realm of approximations. Borel, at the
moment of his probabilistic turn thirty years earlier, expressed
himself similarly when he asserted that a coefficient of probability
constituted the clearest answer to many questions, an answer which
corresponded to an absolutely tangible reality, and when he was ironic
about reluctant spirits who preferred
certainty, and who would perhaps also prefer that 2 plus 2 make 5.\footnote{
  `[\dots ] peut-{\^e}tre aussi que 2 et 2 fissent 5.'}

I refer the reader to the aforementioned studies for precise details
on these questions. What I would like to consider here is how Borel
had combined his research on calculus of probability with the
considerations of his predecessor. From his very first paper, Borel
announced that he adopted the conventionalism of Poincar{\'e}
(\cite{Borel1905}, p.123 ). But his aim was to illustrate the role
that the (then novel) Lebesgue integral and measure theory could play,
after he discovered with amazement their use in \cite{Wiman1900} by
the Swedish mathematician Anders Wiman (on this subject, see
\cite{DurandMazliak2011}).

\begin{quote}`The methods adopted by Mr. Lebesgue allow us to examine
  [\dots ] questions of probability that appear inaccessible to the
  classical methods of integration. Moreover, in the simpler cases, it
  suffices to use the theory of sets I called `measurable',
  and which Mr. Lebesgue later termed `measurable (B)'; the first use of
  this theory of measurable sets for the calculation of probability
  is due, I believe, to Mr. Wiman.'\footnote{ `Les
    m{\'e}thodes de M. Lebesgue permettent d'{\'e}tudier [...] des questions
    de probabilit{\'e}s qui paraissent inaccessibles par les proc{\'e}d{\'e}s
    d'int{\'e}gration classique. D'ailleurs, dans les cas particuliers les
    plus simples, il suffira de se servir de la th{\'e}orie des ensembles
    que j'avais appel{\'e}s mesurables et auxquels M.~Lebesgue a donn{\'e} le
    nom de mesurables (B); l'application de cette th{\'e}orie des
    ensembles mesurables au calcul des probabilit{\'e}s a {\'e}t{\'e}, {\`a} ma
    connaissance, faite pour la premi{\`e}re fois par M. Wiman.'
    (\cite{Borel1905}, p. 126)}\end{quote}

I shall not deal here with the radical transformations the Lebesgue
integral brought to analysis at the beginning of the twentieth
century. For a broad overview, one consult
\cite{Hawkins1975}. Nevertheless, for the sake of completeness, let us
say at least a few words about Borel's role in the elaboration of this
theory.

In his thesis dealing with questions of the extension of analytic
functions, Borel invented a new concept of analytic extension, more
general than that of Weierstrass, using a great deal of geometric
imagination. In the course of his proof, he proved that a countable
subset of an interval can be covered by a sequence of intervals with
total length as small as one wants. This was probably the first
appearance of a $\sigma$-additivity argument for the linear measure of
sets. In subsequent years, Borel considerably fleshed out his
construction, in particular in his work \cite{Borel1898}, by
introducing the notion of measurable set and of measure based on
$\sigma$-additivity. These concepts had however a limited extension
for Borel as he considered only explicit sets obtained by countable
unions and complementary sets, forcing him to make the shaky
suggestion that one should attribute a measure {\it inferior to
  $\alpha$} to any subset of a measurable set with measure $\alpha$.

One had to wait for Lebesgue's thesis and the publication of his Note
\cite{Lebesgue1901}, in which he introduced a new conception for
integration, for the notion of measurable set to reach its full power,
on which is based the remarkable flexibility of the integral exploited
by Borel in his paper \cite{Borel1905}. There he showed in particular
how the use of Lebesgue's integral can allow one to give meaning to
some questions formulated in a probabilistic way.  One of the most
simple of these is that of assigning zero probability to the choice of
a rational number when drawing a real number at random from the
interval [0,1]. Let us insist on the fact that for Borel, it was more
critical that the Lebesgue integral give a question meaning, than to
provide its answer. We see here Borel being absolutely in line with
Poincar{\'e}'s conventionalism, but the choice of the convention
(identifying the probability with the measure of a subset of [0,1]) is
based on mathematical considerations.

However it is above all in his long paper of 1906 on the kinetic
theory of gases (\cite{Borel1906}) that Borel would fit in Poincar{\'e}'s
heritage, at the same time introducing new considerations showing that
he was also striking out on his own. He often insisted on the
difference between his approach and that of Poincar{\'e} (\cite{Borel1906},
p.11, note 2).
 
Borel's aim was to provide a genuine mathematical model for Maxwell's
theory in order to satisfy mathematicians.
 
\begin{quote}`I would like to address all those who shared Bertrand's
  opinion about the kinetic theory of gases, that the problems of
  probability are similar to the problem of finding the captain's age
  when you know the height of the mainmast. If their scruples are
  somewhat justified in virtue of the fact that you cannot fault a
  mathematician for his love of rigor, nevertheless, it does not seem
  to me impossible to give them satisfaction. Such is the aim of the
  following pages: they do not advance the theory in any real sense
  from a physical point of view, but perhaps they will
  convince numerous mathematicians of its interest, and by
  increasing the number of investigators, will contribute
  indirectly to
  its development. Should this happen, these pages will not have been
  useless, independently of the aesthetic interest inherent to any
  logical construction.'\footnote{ `Je voudrais m'adresser {\`a} tous
    ceux qui, au sujet de la th{\'e}orie cin{\'e}tique des gaz, partagent
    l'opinion de Bertrand que les probl{\`e}mes de probabilit{\'e} sont
    semblables au probl{\`e}me de trouver l'{\^a}ge du capitaine quand on
    conna{\^\i}t la hauteur du grand m{\^a}t. Si leurs scrupules sont justifi{\'e}s
    jusqu'{\`a} un certain point parce qu'on ne peut reprocher {\`a} un
    math{\'e}maticien son amour de la rigueur, il ne me semble cependant
    pas impossible de les contenter.'

	`C'est le but des pages qui suivent : elles ne font faire aucun progr{\`e}s r{\'e}el {\`a} la th{\'e}orie du point de vue physique; mais elles arriveront peut {\^e}tre {\`a} convaincre plusieurs math{\'e}maticiens de son int{\'e}r{\^e}t, et, en augmentant le nombre de chercheurs, contribueront indirectement {\`a} son d{\'e}veloppement. Si c'est le cas, elles n'auront pas {\'e}t{\'e} inutiles, ind{\'e}pendamment de l'int{\'e}r{\^e}t esth{\'e}tique pr{\'e}sent dans toute construction logique.' (\cite{Borel1906}, p. 10)} \end{quote}

Thus a motivation for Borel's agenda was that he regarded the various
considerations of Poincar{\'e} on kinetic theory as insufficient
to convince mathematicians of its interest. Let us observe in passing that
Poincar{\'e}, that same year 1906, wrote a new paper for the {\it Journal
  de Physique}, where he studied the notion of entropy in the kinetic
theory of gases (\cite{Poincare1906}); there was probably no direct
link between Poincar{\'e}'s and Borel's publications, which treat different
questions.

Borel began his paper by returning to one of Poincar{\'e}'s major themes:
the distribution of the small planets. However, he approached it
from a new angle (see below). He then applied the results he obtained
to the construction of a mathematical model from which Maxwell's
distribution law can be deduced. Borel's fundamental
idea was that in the phase space where coordinates are the
speeds of $n$ molecules, the sum of squares of these speeds at a
given time $t$ is equal (or, more exactly, proportional) to $n$ times
the mean kinetic energy, so that the point representing the system of
the speeds belongs to a sphere with a radius proportional to $\sqrt
n$. Borel went on to perform an asymptotic study of the uniform measure on
the ball with radius $\sqrt n$ in dimension $n$. I refer the
interested reader to \cite{Vonplato1991}, \cite{Vonplato1994} et
\cite{Mazliak2011} for details, and shall restrict
myself to some comments on the first part of \cite{Borel1906}, dealing
with the small planets.

%
%

Considering a circle on which there are points representing the
longitudinal position of the small planets, Borel posed the
following question: What is the probability for all the small planets
to be situated on the same half-circle fixed in advance? As Borel
noted, if one had perfect knowledge of the positions of the planets,
there would be no reason to invoke probabilities, as one could directly
assert whether the event was realized or not. He argued that it was
necessary to restate the question to give it well-defined
probabilistic meaning according to a selected convention. The
simplest convention would be to assume that the probability for each
planet on the chosen half-circle $C_1$ is equal to the
probability of being on the complementary half-circle (and therefore equal to
1/2), and that the different planets are situated independently with
respect to each other. In this case, naturally, if there are $n$
planets, the desired result is $1/2^n$. However, if this independence
was more or less tacitly considered by Poincar{\'e}, Borel challenged it
as being questionable, the planets having clearly mutual influences
(\cite{Borel1906}, p. 12), and so he sought to forgo this
hypothesis. Progressively enlarging his initial problem, Borel arrived at
the following asymptotic formulation (\cite{Borel1906}, p. 15):

\begin{quote}Problem C. - Given the mean motions of $n$ small
  planets to within $\varepsilon$ and their exact initial positions,
  one denotes by $\overline \omega$ the probability that, at a time
  $t$ chosen at random in an interval $a,b$, every corresponding point $P$
  is in $C_1$. What is the limit to which
  $\overline \omega$ tends when the interval $a,b$ increases
  indefinitely?\footnote{ `Probl{\`e}me C. - Connaissant {\`a} $\varepsilon$
    pr{\`e}s les moyens mouvements des $n$ petites plan{\`e}tes et connaissant
    exactement leurs positions initiales, on d{\'e}signe par $\overline
    \omega$ la probabilit{\'e} pour qu'{\`a} une {\'e}poque $t$ choisie
    arbitrairement dans un intervalle $a,b$ tous les points $P$
    correspondants soient sur $C_1$. Quelle est la limite vers
    laquelle tend $\overline \omega$ lorsque l'intervalle $a,b$
    augmente ind{\'e}finiment?' (\cite{Borel1906}, p. 15)}\end{quote}

Borel could then implement a method of arbitrary functions in
dimension $n$ without supposing the initial independence of the
motions of the planets, and prove that asymptotically the desired
probability was $1/2^n$, a type of ergodic theorem which showed an
asymptotic independence he would also show in the case of his model
for the kinetic theory of gases having Gaussian distributions as
limits. A rather curious detail is that the result given by Borel,
proving the convergence of the uniform distribution on a sphere with
radius $\sqrt n$ in $n$ dimensional space to independent Gaussian
variables, is today called {\it Poincar{\'e}'s Lemma}, even though it is
entirely absent from the works of Poincar{\'e} (for complements about this
strange fact, see \cite{Mazliak2011} and the references therein).

I have not been able to make out clearly whether Poincar{\'e} was ever
interested in the research and the work of his successor in the field
of probability. The only sign which might indicate at least a passing
interest is the fact that he agreed to write an article entitled ``Le
hasard'' (\cite{Poincare1907}) for the {\it Revue du Mois}. But, to my
knowledge, Poincar{\'e} did not comment on Borel's work, and what is even
more surprising, if one remembers Poincar{\'e}'s work at the beginning of
the 1890s, he showed no interest in measure theory as applied to the
mathematics of randomness. Poincar{\'e}, here again, stood poised on the
threshold of a domain he helped to create, but did not enter.

\subsection{Markovian descent}

In order to complete this outline of Poincar{\'e}'s probabilistic
heritage, let us finally consider what may have been the most amazing
consequence of his work: the dazzling development, since the end of
the 1920s, of the theory of Markov chains and Markov processes. This
story has already been set out in several texts and I shall again
restrict myself to comments on only the most salient points, referring
the reader elsewhere for more information.

We have already evoked Poincar{\'e}'s investigations of card shuffling and
the fact that in his proof of the convergence to the uniform
distribution in \cite{Poincare1912}, he used an algebraic method with
limited exploitation of the probabilistic structure of the
model. Borel, an attentive reader, immediately realized this and wrote
a note, asking Poincar{\'e} to present it to the Academy for the {\it
  Comptes-Rendus de l'Acad{\'e}mie des Sciences} (the only letter from
Borel placed online on the website of the Archives Poincar{\'e}
\footnote{http://www.univ-nancy2.fr/poincare/chp/ }).

Borel wrote to his colleague on 29 December 1911:

\begin{quote}`I have just read the book you kindly sent me; I do
  not need to tell you how much the new sections interested me, in
  particular your theory of card shuffling. I tried to render it
  accessible to those unfamiliar with complex numbers, and
  it seems to me that in doing so I reached a slightly more general
  proposition. If it is new, and if you find it interesting, I would
  ask you to communicate the attached note.\footnote{
    `Je viens de lire le livre que vous avez eu l'amabilit{\'e} de me
    faire envoyer; je n'ai pas besoin de vous dire combien les parties
    neuves m'ont int{\'e}ress{\'e}, en particulier votre th{\'e}orie du
    battage. J'ai essay{\'e} de la mettre {\`a} la port{\'e}e de ceux qui ne sont
    pas familiers avec les nombres complexes, et il m'a sembl{\'e} que
    j'obtenais ainsi une proposition un peu plus g{\'e}n{\'e}rale. Si elle est
    nouvelle, et si elle vous para{\^\i}t int{\'e}ressante, je vous demanderai
    de communiquer la note ci-jointe.'}\end{quote}

Poincar{\'e} acted immediately and the note was presented on 2 January
1912. Borel's method in \cite{Borel1912} was in fact an extension of
the elementary one used by Poincar{\'e} to treat the case of two cards,
where one looks at the evolution of the successive means in the course
of time. This method later became the standard in proofs of the exponential
convergence of an irreducible finite Markov chain towards its
stationary distribution (see for instance \cite{Billingsley1995},
p.131). Here the stationary distribution is uniform due to the
reversible character of the chain. Borel even gave himself the
satisfaction of introducing a dependence on time (the chain becoming
inhomogeneous).

He considered the regular case where there exists an $\varepsilon$
such that, at every moment, the transition probabilities of one
permutation to another at a subsequent time are all greater than some
$\varepsilon$.  In Borel's notation, let $p_{j,n}$ be the probability
of the $j$-th possible permutation of the cards before the $n$-th
operation. Denoting by $\alpha_{j,h,n}$ the probability for $A_h$ to
be replaced by $A_j$ during the $n$-th operation, one
has $$p_{j,n+1}=\sum_{h=1}^{h=k}\alpha_{j,h,n}p_{h,n}$$ with the
constraint $\sum_{h=1}^k\alpha_{j,h,n}=1$ where $k$ denotes the number
of possible permutations. Let us immediately observe that $P_n$ and
$p_n$, the largest and the smallest of the $p_{j,n}$, form two
sequences, respectively nonincreasing and nondecreasing. Let $P$ and
$p$ denote their limits. For a given $\eta >0$, one may choose $n$ for
which $P_n\le P+\eta $, and therefore the $p_{j,n}$ are inferior to
$P+\eta$. After $N$ operations one can write
$$p_{j,n+N}=\sum_{h=1}^{h=k}\beta_{j,h,n}p_{h,n}\hskip 5pt ,\hskip 5pt \sum_{h=1}^{h=k}\beta_{j,h,n}=1$$ where the $\beta_.$  are the transition probabilities between time $n$ and time $n+N$, each being greater than $\varepsilon$ by hypothesis. 

Let us consider the smallest of the $p_{h,n}$, $p_{h_0,n}$ so that
$p_n=p_{h_0,n}\le p$. For the sake of simplicity, let us denote by
$\beta$ its coefficient $\beta_{j,h_0,n}$; by hypothesis, $\beta \ge
\varepsilon$. Let us observe that $\sum_{h=1, h\neq
  h_0}^{h=k}\beta_{j,h,n}= 1-\beta$. Therefore, one can write, by
choosing $j$ such that $p_{j,n+N}$ is superior or equal to $P$,
$$
P\le p_{j,n+N}\le \beta  p + (1-\beta )(P+\eta )= P+(1-\beta )\eta -\beta (P-p)
$$ 
and hence 
$$
P-p \le \frac{1-\beta}{\beta}\eta \le \frac{1-\varepsilon}{\varepsilon}\eta .
$$ 
$\eta$ being arbitrary small, one concludes that $P=p$ and therefore
that asymptotically the $p_{j,n}$ become all equal to $1/k$. Let us
observe in passing that in these blessed years when it was permitted
to publish mistakes, Borel erred in writing his inequality,
considering $\varepsilon$ instead of the number we called $\beta$, a
fact which naturally did not change anything in the final result.

Nobody seemed to have paid attention to Borel's note: when these
results were rediscovered by L{\'e}vy and then Hadamard in the 1920s,
neither of them had the slightest idea of its existence (on this
subject see the letters from L{\'e}vy to Fr{\'e}chet \cite{Barbut2004}, pp.137
to 141).

We must next skip five years and cover several hundreds of kilometers
east in order to see a new protagonist coming on stage, the Czech
mathematician Bohuslav {\sc Hostinsk\' y}. Moreover, as if it were not
enough that we must invoke an unknown mathematician, we must first say
a few words about an unknown philosopher. Indeed the man who may have
been, together with Borel, the most attentive contemporary reader
of Poincar{\'e}'s texts on probability was another Czech, the philosopher
Karel {\sc Vorovka} (1879-1929), whose influence on Hostinsk\'y was
decisive.

It is not possible here to discuss this singular figure in detail and
I shall therefore restrict myself to giving some elements explaining
how he had got involved in this melting pot. An interesting and very
complete study on Vorovka was published in Czech some years ago
\cite{Pavlincova2010} and hopefully it will become more accessible in
a more widely known language. Some complements can also be found in
\cite{Mazliak2007b} and the references therein. Two reasons explain
this general ignorance of Vorovka: the fact that his works, mostly in
Czech, were never translated, and also that his early death
precluded the collection of his ideas in a large-scale work. Placing
himself in the tradition of Bernhard Bolzano (1771-1848), the major
figure of the philosophical scene in Prague during the 19th century,
Vorovka sought an approach combining both his
strong mathematical education and a rather strict religious
philosophy, an original syncretism of empiricism and idealism which
had close links with the thought of the hero of the
Czechoslovak independence, T.G.~Masaryk, and with American
pragmatist philosophy, in which he had much been interested.

Vorovka's discovery of Poincar{\'e}'s philosophical writings at the beginning of
the 20th century was a real revelation: he drew from them
the conviction that the scientific discoveries
at the end of the 19th century, especially in physics, compelled a
reconsideration of the question of free will. Vorovka
showed a real originality in that he did not content himself with
principles, but closely studied the mathematical problems raised by
the theory of probability. He was a diligent reader of Bertrand's
textbook, of Borel's texts, but also of Markov's works, publishing
several works inspired by papers of the Russian mathematician (see
\cite{Markov1905}, \cite{Vorovka1912a}, \cite{Vorovka1912b},
\cite{Vorovka1914}). At the time when he was granted tenure at the
Czech University in Prague, around the year 1910, Vorovka met the
mathematician Bohuslav Hostinsk\'y, who had just returned to Bohemia
after a period of research in Paris. In Hostinsk\'y's own words (see
\cite{Hostinsky1929}), it is through the discussions he had with
Vorovka that he learned about Poincar{\'e}'s works, and he began to reflect
upon the calculus of probability, a domain somewhat remote from his
original field of research (differential geometry).

Following Ji\v{r}{\'\i}\hskip 2pt Ber{\'a}nek, who was one of the last
assistants to Hostinsk\'y after the Second World War at the University
of Brno, another source of the latter's interest in the calculus of
probability is found in his reading of the 1911 paper by Paul and
Tanya Ehrenfest on Statistical Mechanics for the {\it Encyklop{\"a}die der
  Mathematischen Wissenschaften}, translated and completed by Borel
for the French version of the {\it Encyclop{\'e}die des Sciences
  Math{\'e}matiques} \cite{Ehrenfest1915}.

Ber{\'a}nek wrote (\cite{Beranek1984}) that this paper, whose impact was considerable,
\begin{quote} `put the emphasis on statistical methods in physics,
  along with geometrical methods, mainly in connection with the works
  of L. Boltzmann on kinetic gas theory. Boltzmann's work sustained
  discussions and controversies concerning the correctness and
  legitimacy of his mathematical methods.
Hostinsk\'y, as he mentioned
  later, began to study Boltzmann's works in 1915, and to take an
  interest in the efforts made to provide precise mathematical bases
  for the kinetic theory. The central point of these efforts implied a reexamination of some fundamental questions of the theory of
  probability. Hostinsk\'y was especially impressed by the fundamental
  works of H. Poincar{\'e} on the bases of probability calculus which
  opened the way for new methods necessary for the improvement of
  kinetic theory. For this reason, around 1917, Hostinsk\'y began to
 study in earnest questions in the calculus of probability\dots
  '\footnote{ `[\dots ] mettait l'accent sur les m{\'e}thodes statistiques
    en physique, {\`a} c{\^o}t{\'e} des m{\'e}thodes g{\'e}om{\'e}triques, principalement en
    relation avec les travaux de L. Boltzmann sur la th{\'e}orie cin{\'e}tique
    des gaz. Sur ceux-ci furent men{\'e}es discussions et controverses, au
    sujet de l'exactitude et de la l{\'e}gitimit{\'e} des m{\'e}thodes
    math{\'e}matiques employ{\'e}es. Hostinsk\'y, comme il l'a lui m{\^e}me
    mentionn{\'e}, commen\c ca {\`a} partir de 1915 {\`a} {\'e}tudier les travaux de
    Boltzmann et {\`a} s'int{\'e}resser aux efforts qui {\'e}taient faits pour
    donner {\`a} la th{\'e}orie cin{\'e}tique des bases math{\'e}matiques pr{\'e}cises. Le
    point central de ceux-ci n{\'e}cessitait un nouvel examen de certaines
    questions fondamentales de la th{\'e}orie des
    probabilit{\'e}s. Hostinsk\'y fut particuli{\`e}rement impressionn{\'e} {\`a} ce
    sujet par les travaux fondamentaux de H. Poincar{\'e} sur les
    fondements du calcul des probabilit{\'e}s qui ouvraient la voie {\`a} de
    nouvelles m{\'e}thodes n{\'e}cessaires pour le perfectionnement de la
    th{\'e}orie cin{\'e}tique. Pour cette raison, vers 1917, Hostinsk\'y
    commen\c ca {\`a} s'occuper s{\'e}rieusement de questions de calcul des
    probabilit{\'e}s\dots '} \end{quote}

The fact that Hostinsk\'y began to deal seriously with probability in
1917 is attested by his own diary, kept in the archives of Masaryk's
university in Brno. The diary's entries concern only comments on differential geometry until 1917. On 10 January 1917,
Hostinsk\'y made some observations on the study of card shuffling by
Poincar{\'e}, following \cite{Poincare1912}, and on January 18th he took up problems
of lottery. A first paper appeared some months later in the {\it
  Rozpravy \v{C}esk{\'e} Akademie} dealing with the problem of Buffon's
needle \cite{Hostinsky1917}.

The problem of Buffon's needle is a classic of the calculus of
probability and Hostinsk\'y began by expounding it:

\begin{quote} `A cylindrical needle is thrown on a horizontal floor,
  on which are traced equidistant parallels; the distance $2a$ between
  two successive parallels is supposed larger than the length $2b$ of
  the needle. What is the probability that the needle meet one of the
  parallels?'\footnote{ `On lance une aiguille cylindrique sur un plan
    horizontal, o{\`u} sont trac{\'e}es des parall{\`e}les {\'e}quidistantes; la
    distance $2a$ de deux parall{\`e}les voisines est suppos{\'e}e plus grande
    que la longueur $2b$ de l'aiguille. Quelle est la probabilit{\'e} pour
    que l'aiguille rencontre l'une des parall{\`e}les?'} \end{quote}
 
Buffon had proposed a solution whose numerical result $\frac{2b}{\pi
  a}$, in which $\pi$ was present, was a source of numerous
propositions for an `experimental' calculation of $\pi$. But, in fact,
Buffon's proof was based on the hypothesis that the needle center
could be located anywhere on the plane, and
Hostinsk\'y, in a second critical part of his paper, mentioned the
dubious nature of such an hypothesis, just as Carvallo had done before
him in 1912. An experimental device could only take the form of a
table of limited size, and it is then clear that, depending on the choice
of a small square $C_1$ at the center of the table or another square
$C_2$ on the edge of the table with the same area, the probability
$p_1$ that the center of the needle belongs to $C_1$ and the
probability $p_2$ that it belongs to $C_2$ cannot be the same: indeed,
$C_2$ is strongly subject to the constraint that the needle does not
fall from the table, but $C_1$ is very weakly so constrained, so that intuitively
one should have $p_1>>p_2$.
 
Hostinsk\'y therefore considered it indispensable to suppose unknown
the {\it a priori} distribution of the localization of the needle. It
is an {\it unknown} distribution (with density) $f(x,y)dxdy$. But,
mentioned Hostinsk\'y, Poincar{\'e} also, in the resolution of several
problems of probability, allowed the use of such an {\it arbitrary}
density and observed that in some situations this function would not
be present in the final result. Hostinsk\'y proposed to prove that if
a domain $A$ of the space is segmented in $m$ elementary domains with
the same volume $\varepsilon$, and containing each a white part with
volume $\lambda \varepsilon$ and a black part with volume $(1-\lambda
) \varepsilon$ (where $0<\lambda <1$), then for any sufficiently
regular function $\varphi (x,y,z)$, the integral on the white parts
will asymptotically (when $m$ tends to infinity) be equal to $\lambda$
times the integral of $\varphi$ on $A$.
 
Hostinsk\'y then applied this result in order to propose a {\it new}
solution to the problem of the needle. Instead of Buffon's unrealistic
hypothesis, he supposed that the center of the needle is compelled to
fall in a square with side $2na$, $n\in \mathbb N$, with a density of
probability given by an unknown function $\varphi$ (which he supposed
to have bounded derivatives) and kept on the other hand the second
hypothesis concerning the uniform distribution of the angle $\omega$
of the needle with respect to the parallels. This being set out,
dividing the domain of integration $0<x<2na, 0<y<2na, 0<\omega
<\frac{\pi }{2}$ in $n^2$ subdomains (by partitioning the values of
$x$ and $y$ with respect to the multiples of $a$), each small domain
is itself divided into two parts (corresponding to the fact that the
needle intersects [{\it white} part] or does not intersect [{\it
  black} part] the corresponding parallel). The ratio of their
respective volumes to the total volume of the subdomain is constant
and equal, for the white part, to $\frac{2b}{\pi a}$. An application
of the previous theorem then allows one to assert that we obtain the
desired probability, at least asymptotically when $n$ tends towards
infinity.
 
In the Spring of 1920, seeking to benefit from the sympathy of French
public opinion towards the young Czechoslovakia, Hostinsk\'y had sent
to {\'E}mile Picard the translation of his paper and Picard proposed
immediately (18 April 1920) to include it in the M{\'e}langes of the {\it
  Bulletin des Sciences Math{\'e}matiques.} This slightly revised version
of the paper of 1917 was published at the end of 1920 and Maurice
Fr{\'e}chet, who had just arrived in Strasbourg and considered himself as
a missionary \cite{SiegmundSchultze2005} read it with attention, as he
mentioned in a subsequent letter to Hostinsk\'y, dated 7 November
1920, in which he congratulated him on having obtained such a {\it
  positive result}.

As we have just explained, following Poincar{\'e}'s example, Hostinsk\'y
required that the function $\varphi$ admit a uniformly bounded
derivative in the domain $A$ in order to obtain an upper bound for the
difference between the maximum and the minimum of $\varphi$ on each of
the small domains. But Fr{\'e}chet, when he read the paper, rightly
realized that as only an estimation of the integrals of $\varphi$ on
these domains was needed, the simultaneous convergence of the superior
and inferior Darboux sums towards the integral of $\varphi$ allowed
one to obtain the desired result with $\varphi$
Riemann-integrable. This is what he wrote, together with the proof, to
Hostinsk\'y on 7 November 1920.

It seems that the former letter refers to
Fr{\'e}chet's initial research on probabilistic questions. It was subsequently published in a
short note in 1921 (\cite{Frechet1921}). Hostinsk\'y answered on 22
December 1920, agreeing with Fr{\'e}chet that the integration hypothesis
was sufficient. He also mentioned that Borel had already suggested that
Poincar{\'e}'s hypothesis could be weakened, supposing only the function
to be continuous. In his textbook on probability published in 1909 by
Hermann \cite{Borel1909}, in which Borel devoted the whole of Chapter
VIII to the introduction of arbitrary functions by considering both
Poincar{\'e}'s examples of the roulette wheel and of the small planets on
the Zodiac, Borel noted that the hypothesis of continuity was
sufficient to apply Poincar{\'e}'s method. Fr{\'e}chet was to include
Hostinsk\'y's observation in his note in 1921 \cite{Frechet1921}
(where he emphasized that it had been inspired by the latter after
having read his paper on Buffon's needle). In \cite{Frechet1921}, he
mentioned Borel's work to emphasize that hypotheses of
both continuity (Borel) and derivability (Poincar{\'e}) were useless and that
Riemann-integrability was sufficient.

His friendly relationship with Fr{\'e}chet encouraged Hostinsk\'y to continue
his probabilistic studies, and this leads us eventually to the last
step of this long journey, introducing Jacques Hadamard
(1865-1963). The presence of this name in our story may seem quite
strange, and, in fact, Hadamard was interested in probabilities only
during one semester of the academic year 1927-1928. He had never
considered them before, and would never do so again, showing even some
irritation towards L{\'e}vy, one of the disciples of whom he was most fond,
when he `wasted' his mathematical talent in the 1920s and left the
royal path of functional analysis for the calculus of probability.

Following Poincar{\'e}'s example, Hadamard always kept in mind
physical theories from which he intended to extract new mathematical
problems. It was from this perspective that he had written the
aforementioned review of Gibbs' book in 1906.

When Hadamard began writing up his course of
analysis at the {\'E}cole Polytechnique in the 1920s (published by Hermann in two
volumes in 1926 and 1930), he had to prepare some lectures in
1927 on
probability theory, and he took up Poincar{\'e}'s example of card
shuffling. On this occasion, he recovered Borel's method of successive
means and published in 1927 a note in the {\it Comptes-Rendus de
  l'Acad{\'e}mie des sciences de Paris} \cite{Hadamard1927}. Soon after
that, Hostinsk\'y discovered Hadamard's note and sent an extension of
it for publication in the {\it Comptes-Rendus}
\cite{Hostinsky1928}, which appeared in the first weeks of
1928. There, for the first time, before everyone {- except for
  Bachelier, but, alas, who had ever read Bachelier among the
  mathematicians I write about! - } and especially before Kolmogorov,
Hostinsk\'y introduced a Markovian model in continuous time. {At the
  Bologna congress in September 1928, it was realized that Poincar{\'e}'s
  card-shuffling studies were in fact a special case of the model
  of variables in chain introduced by Markov in 1906, and developed in
  several of his posterior publications as well as by Bernstein in
  1926, but which were largely ignored outside Russia. }The attention
this drew, in particular at the Congress in Bologna in 1928,
inaugurated intense activity on these questions which continued
throughout the 1930s, a story brilliantly recounted in \cite{Bru2003}
to which I refer the interested reader. This unexpected crowning of
Poincar{\'e}'s efforts seems to be a perfect moment to take leave of the
master.

\section*{Conclusion}

Poincar{\'e} lived during that very specific moment in the history of
science when randomness, in a more and more insistent way, challenged
the beautiful deterministic edifice of Newton's and Laplace's cosmology
which had dominated scientific thinking for centuries. A conference by
Paul Langevin in 1913 \cite{Langevin1913} shows the extent of this
challenge, paralleling the introduction of probabilities and a drastic
change in our comprehension of the structural laws of matter. Such a
penetrating mind as Poincar{\'e} could not have lived this irruption
otherwise than as a traumatic one, that he had to face with the means
he had at his disposal. These means, as we saw, had not yet reached
the degree of power necessary to deal with many problems raised by
modern physics. Let us recall one of the master's apothegms:

\begin{quote} `Physics gives us not only an occasion for problem
  solving: it helps us find the means of solution, and in two ways. It
  points us in the direction of the solution, and suggests how to
  reason our way there.'\footnote{ `La physique ne nous donne pas
    seulement l'occasion de r{\'e}soudre des probl{\`e}mes; elle nous aide {\`a}
    en trouver les moyens, et cela de deux mani{\`e}res. Elle nous fait
    pressentir la solution; elle nous sugg{\`e}re des
    raisonnements.'(\cite{Poincare1905}, p. 152)}
\end{quote}

And the new physical theories with which Poincar{\'e} was confronted
suggested developing the theory of probability in the first place - a
suggestion which can be also found, but in a slightly different
perspective, among the problems Hilbert expounded during the Paris
Congress of 1900. Therein lies the apparent paradox which puzzled the
mathematician at the turn of the century: the hesitation and
reluctance in the face of problems raised by statistical mechanics,
the somewhat uncertain attempts to give solid bases to the theory of
probability, the seemingly limited taste for new mathematical
techniques, in particular measure theory and Lebesgue's integration,
though they could have provided decisive tools to tackle numerous
problems. Poincar{\'e}, as we said, remained a man of the 19th century,
maybe in the same way as Klein had mischievously presented Gauss as a
scientist of the 18th century. Naturally, in Gauss's case, the irony
came from the fact that he had lived two thirds of his life in the
19th century, whereas death surprised Poincar{\'e} at the beginning of the
20th century. But we may speculate - although {\it not} here! - on the
manner in which our hero would have adapted to transformations in the
scientific world picture. We have seen that, following the example of
his glorious predecessor, Poincar{\'e} sowed widely, and the spectacular
blossoms of many of his ideas inspired countless researchers after his
death. As for probabilities, I think one can sum up the measure of his
influence as follows: he began to extract the domain from the grey
zone to which it had been confined by almost all French
mathematicians, he initiated methods that flourished when they
integrated more powerful mathematical theories, he convinced
Borel of the importance of certain questions, to the study of which he
eventually devoted an enormous amount of energy. For a rather marginal
subject in Poincar{\'e}'s works, such a contribution appears far from
negligible.

\end{document}